\def\no{\if01}
\def\iftwelvept{\no}

\def\ifusepdf{\no}
\def\ifpsfont{\no}

\documentclass[leqno,12pt]{amsart}
\usepackage{amssymb}
\usepackage{amscd}
\usepackage{latexsym}
\usepackage{verbatim}
\usepackage{mathrsfs}
\usepackage[all]{xy}

\ifusepdf
\usepackage{hyperref}
\else\fi
\ifpsfont
\usepackage[T1]{fontenc}
\usepackage{times}
\else\fi

\iftwelvept
\else\fi

\newtheorem{Theorem}{Theorem}[section]
\newtheorem{Proposition}[Theorem]{Proposition}
\newtheorem{Lemma}[Theorem]{Lemma}
\newtheorem{Corollary}[Theorem]{Corollary}

\newtheorem{Remark}[Theorem]{Remark}

\newcommand{\ZZ}{{\mathbb{Z}}}

\newcommand{\RR}{{\mathbb{R}}}

\newcommand{\SSSS}{\mathscr{S}}
\newcommand{\DDD}{{\mathcal{D}}}
\newcommand{\EEE}{{\mathcal{E}}}

\newcommand{\CCC}{{\mathcal{C}}}

\newcommand{\OO}{{\mathcal{O}}}
\newcommand{\Fun}{\operatorname{Fun}}
\newcommand{\Map}{\operatorname{Map}}

\newcommand{\KKK}{\mathcal{K}}

\newcommand{\lcatinf}{\widehat{\mathcal{C}\textup{at}}_{\infty}}
\newcommand{\UU}{\mathcal{U}}
\newcommand{\DP}{\mathcal{D}_{\operatorname{perf}}}

\newcommand{\DQ}{\mathcal{D}_{\operatorname{qcoh}}}

\newcommand{\DV}{\mathcal{D}_{\operatorname{vect}}}
\newcommand{\Catinf}{\mathcal{C}\textup{at}_{\infty}}

\newcommand{\Hom}{\operatorname{Hom}}
\newcommand{\Ext}{\operatorname{Ext}}

\newcommand{\NNNN}{\operatorname{N}}

\newcommand{\Spec}{\operatorname{Spec}}

\newcommand{\FFF}{\mathbb{F}}

\newcommand{\Catdel}{\textup{Cat}_{\Delta}}
\newcommand{\Setdel}{\textup{Set}_{\Delta}}

\newcommand{\XX}{\mathcal{X}}

\newcommand{\catsmon}{\widehat{\mathcal{C}\textup{at}}_{\infty}^{\textup{sMon}}}
\newcommand{\Proof}{{\sl Proof.}\quad}
\newcommand{\QED}{{\unskip\nobreak\hfil\penalty50\quad\null\nobreak\hfil
{$\Box$}\parfillskip0pt\finalhyphendemerits0\par\medskip}}

\DeclareMathOperator{\Coch}{C^\bullet}

\DeclareMathOperator{\id}{id}
\DeclareMathOperator{\chain}{C_\bullet}

\newcommand{\lotimes}{\otimes^{\bf{L}}}

\newcommand{\stk}[1]{\mathcal {#1}}

\newcommand{\dqc}{\mathcal D_{\mathrm{qcoh}}}
\newcommand{\dpf}{\mathcal D_{\mathrm{perf}}}

\newcommand{\dvb}{\mathcal D_{\mathrm{vect}}}
\newcommand{\tdqc}{\mathrm D_{\mathrm{qcoh}}}
\newcommand{\tdpf}{\mathrm D_{\mathrm{perf}}}

\newcommand{\hhom}{\mathrm{Hom}}

\newcommand{\cpx}{}

\newcommand{\ps}[1]{\mathrm{pr}_{#1*}}

\newcommand{\pl}[1]{\mathrm{pr}_{#1}^*}
\newcommand{\ol}[1]{\overline{#1}}

\newcommand{\bu}{\mbox{1}\hspace{-0.25em}\mbox{l}}
\newcommand{\tH}{\mathrm H}
\newcommand{\mfn}{\Phi}

\newcommand{\co}{\widehat{\OO}}

\DeclareSymbolFont{cmletters}{OML}{cmm}{m}{it}
\DeclareSymbolFont{cmsymbols}{OMS}{cmsy}{m}{n}
\DeclareSymbolFont{cmlargesymbols}{OMX}{cmex}{m}{n}
\DeclareMathSymbol{\myjmath}{\mathord}{cmletters}{"7C}
\DeclareMathSymbol{\myamalg}{\mathbin}{cmsymbols}{"71}
\DeclareMathSymbol{\mycoprod}{\mathop}{cmlargesymbols}{"60}
\let\jmath\myjmath

\let\coprod\mycoprod

\begin{document}

\title[Monoidal $\infty$-category of complexes from Tannakian viewpoint]{Monoidal Infinity Category of complexes \\ from Tannakian viewpoint}



\author{Hiroshi Fukuyama and Isamu Iwanari}
\address{Department of Mathematics, Faculty of Science,
Kyoto University, Kyoto, 606-8502, Japan}
\email{fukuyama@math.kyoto-u.ac.jp}
\address{Mathematical Institute, Tohoku University, Sendai, Miyagi, 980-8578, Japan}
\email{iwanari@math.tohoku.ac.jp}

\thanks{We take this opportunity to correct a mistake made by the publisher, Springer; in the published version the received date is assigned to 10 May 2012,
though we submited this article on 10 May 2010.}


%

\maketitle

\begin{abstract}
In this paper we prove that a morphism between schemes or stacks
naturally corresponds to
a symmetric monoidal functor between stable $\infty$-categories
of quasi-coherent complexes.
It can be viewed as a derived analogue of Tannaka duality.
As a consequence,
we deduce that an algebraic stack satisfying a certain condition
can be recovered from the symmetric monoidal stable $\infty$-category of quasi-coherent complexes with tensor operation.

\end{abstract}

\section{Introduction}\label{Intro}

Let $X$ be a reasonably nice scheme (for example, a noetherian scheme). Many important invariants of $X$
come from the category of complexes of coherent sheaves
on $X$.
Practically speaking, by the category of complexes
we mean the derived category of coherent complexes.
The triangulated category equips some natural additional
structures
such as
tensor structure arising from the derived tensor product, the canonical
$t$-structure, etc.
The symmetric monoidal (tensor) structure naturally
determines
the intersection products on algebraic $K$-theory groups,
and thus this product yields the ring structures
on various cohomology theories.
Let us consider the tensor triangulated category
$(\textup{D}_{\textup{perf}}(X),\lotimes)$
of
perfect complexes on $X$, endowed with the derived tensor product $\lotimes$.
In the remarkable papers \cite{Bal}, \cite{Bal2},
Balmer proved that the tensor
triangulated category $(\textup{D}_{\textup{perf}}(X),\lotimes)$
remembers the scheme $X$, that is to say,
the whole scheme $X$ can be recovered
from tensor triangulated category  
$(\textup{D}_{\textup{perf}}(X),\lotimes)$.
Balmer's  reconstruction uses
the classification of tensor thick subcategories of $\textup{D}_{\textup{perf}}(X)$,
which has been studied by Hopkins \cite{Hop}, Neeman \cite{Nee} and Thomason \cite{Tho}.
Roughly speaking, the reconstruction proceeds as follows.
To the tensor triangulated category
$(\textup{D}_{\textup{perf}}(X),\lotimes)$
he associates a ringed topological space which we shall denote
by $\operatorname{Spec}(\textup{D}_{\textup{perf}}(X),\lotimes)$.
A point on the topological space $\operatorname{Spec}(\textup{D}_{\textup{perf}}(X),\lotimes)$
corresponds to a tensor thick subcategory of $(\textup{D}_{\textup{perf}}(X),\lotimes)$
which satisfies a certain condition.
Making use of Thomason's classification of
 tensor thick subcategories of $(\textup{D}_{\textup{perf}}(X),\lotimes)$
 in terms of algebraic cycles,
 Balmer showed that the ringed space $\operatorname{Spec}(\textup{D}_{\textup{perf}}(X),\lotimes)$ is isomorphic to the ringed space $X$.

We are motivated by
this reconstruction problem and the classical Tannaka
duality.
Our principal idea is to view
the reconstruction of a scheme from the tensor triangulated category of perfect complexes as
a derived analogue of Tannaka duality.
Let $G$ be an affine group scheme over a field $k$.
Then Tannaka duality states that $G$ can be reconstructed
from the tensor abelian category of finite dimensional representations
of $G$.
More precisely, if $S$ is an affine $k$-scheme and
$\Hom_k (S,BG)$ is the groupoid
of $S$-valued points (over $k$), then Tannaka duality gives an equivalence
\[
\Hom_k(S,BG) \longrightarrow \Fun_k (\textup{VB}^{\otimes}(BG),\textup{VB}^{\otimes}(S)):\ f\mapsto f^*
\]
where $BG$ is the classifying stack of $G$, $\textup{VB}^{\otimes}(\bullet)$
denotes the tensor exact category of vector bundles,
and $\Fun_k(\textup{VB}^{\otimes}(BG),\textup{VB}^{\otimes}(S))$
is the category of tensor exact $k$-linear functors.
(See \cite{DM}, \cite{Saa} for the precise statement.) 

Let $X$ be a scheme or algebraic stack
(satisfying a certain condition) and let $\textup{D}_{\textup{perf}}(X)$ be
the triangulated category of perfect complexes on $X$.
Our main goal is (roughly speaking) to establish
a derived analogue of Tannaka duality (Theorem~\ref{supermain}),
which relates
the category of morphisms $S\to X$ from a scheme $S$ to an algebraic stack
$X$
with the category of exact functors $\textup{D}_{\textup{qcoh}}(X)\to \textup{D}_{\textup{qcoh}}(S)$ that preserves derived tensor products.
Besides the appealing Tannakian viewpoint
our approach has the virtue of recovering
rich data.
Thick subcategories and tensor thick subcategories
of a (tensor) triangulated category give rise to
localizations.
For a (nice) scheme $X$,
localizations of the triangulated category $\textup{D}_{\textup{perf}}(X)$
arising from Zariski open sets are described in terms
of tensor thick subcategories, and it enables one to reconstruct $X$.
However,
if a (tensor)
triangulated category $D$ is the derived category
arising from algebraic stacks (including the derived category
of complexes of representations of an algebraic group)
and representations of quivers,
the data of tensor thick subcategories in  $D$
is not enough to recover the original sources such as stacks
and quivers, and they happen to be trivial.
For instance, if $X$ is a Deligne-Mumford stack
satisfying a certain condition,
the recent result of Krishna \cite{Kri} shows that
only the coarse moduli space $M$ for $X$ can be recovered from
the data of tensor thick subcategories in $\textup{D}_{\textup{perf}}(X)$.
In our Tannakian approach,
we treat the data arising from symmetric monoidal functors
which are not necessarily localizations.
As a consequence,
our reconstruction is applicable to a fairly large class of 
Deligne-Mumford stacks.
The stabilizer group at a point on a stack is described as the
automorphism group of monoidal natural transformations.

One noteworthy feature of our approach is the usage of
higher category theory.
The natural machinery of higher categories
allows us
to formulate and study our derived Tannaka duality (Theorem~\ref{supermain}).
In addition, it enables us to prove a categorical characterization
of derived functors associated to 
morphisms of schemes and stacks (Theorem~\ref{geomono}).
In particular, the characterization describes the clear relationship
between the automorphism group of a variety (or stack) $X$
and the group of autoequivalences of  the ``derived category'' of $X$.
In order to treat symmetric monoidal functors
and realize the derived Tannaka formalism
we shall replace the triangulated category $\textup{D}_{\textup{perf}}(X)$
by ``enhanced'' (higher) category $\mathcal{D}_{\textup{perf}}(X)$.
There are some candidates
which provide the frameworks dealing with such enhanced higher categories:
triangulated derivators, dg-categories, stable simplicial categories,
stable Segal categories and stable $\infty$-categories (quasi-categories), etc.
We use the theory of $\infty$-categories (quasi-categories)
which has been extensively developed by Joyal and Lurie \cite{J}, \cite{HTT}. In addition, many parts of this paper depend on
the theory of $\infty$-categories and theorems such as derived Morita theory
\cite{To}, \cite{BFN}
build on the higher category theory.

There should be various viewpoints and possible formulations for
realizing tannakian
ideas and phenomena in higher category theory.
A tannakian idea of higher category theory appears in \cite{TH}.
Also, we would like to invite the reader's attention 
to the recent works \cite{Tan}, \cite{DAG8}, \cite{Wall}.

The contents of this paper are organized as follows.
In Section \ref{sec:2},  we begin by reviewing the basic notions
on $\infty$-categories in the sense of \cite{J} and \cite{HTT}
and we give preliminaries related to our study.
Section \ref{sec:3}, \ref{sec:4}, \ref{sec:5} contain our main discussions.
In Section \ref{sec:3},
we prove some lemmas concerning Kan extensions
for functors between $\infty$-categories.
This is generalized to a version of symmetric monoidal
functors in Section \ref{sec:5} before the proof of main results.
In Section \ref{sec:4}, we then 
proceed to prove some key property of symmetric monoidal functors
between stable symmetric monoidal
$\infty$-categories of quasi-coherent complexes.
We prove that a colimit-preserving functor
$\mathcal{D}_{\textup{qcoh}}(X) \to \mathcal{D}_{\textup{qcoh}}(Y)$
which commutes with symmetric monoidal operation,
preserves vector bundles.
This property is vaguely reminiscent of
``tannakian phenomenon'' which relates schemes (stacks)
with stable $\infty$-categories of quasi-coherent complexes over them.
In our study derived Morita theory plays an important role.
In Section \ref{sec:5}, 
applying results of Section \ref{sec:3} and \ref{sec:4}, we will prove
main results of this paper.

%
%
%
%
%
%
%
%
%
%
%
%
%
%
%
%

\section{Preliminaries and $\infty$-category of complexes}\label{sec:2}

In this section, we will fix notion and convention and
prepare the settings.
We begin by reviewing the theory of $\infty$-categories which
we will use in the course of this paper.
Roughly speaking, an $(\infty,1)$-category
or simply an $\infty$-category is a weak $\infty$-category
whose $n$-morphisms are invertible for $n>1$.
At present, there are at least four approaches to such a theory:
simplicial categories, Segal categories, complete Segal spaces
and quasi-categories. It is known that all four theories are equivalent.
In other words, each theory is linked to one another via a
Quillen equivalence (see \cite{JT}, \cite{Berg}).
Among them, we use the theory of quasi-categories (\cite{JJ}, \cite{J}, \cite{HTT}), which we shall call $\infty$-categories.
We review basic definitions and facts on quasi-categories
for the convenience of the reader.
However, it is an almost impossible task to present a rapid
overview of all materials
\cite{JJ}, \cite{J}, \cite{HTT}, \cite{HA}
and thus our review is a quick introduction to
basic notions on quasi-categories, which are appearing in
the first Chapter of \cite{HTT}.
Therefore we refer to the books
\cite{HTT}, \cite{HA}
as the general reference of the theory of quasi-categories.

\subsection{$\infty$-categories}\label{subsec:2-1}
Let us recall the definition of a quasi-category.
A (small) {\it quasi-category} $S$ is a (small) simplicial set
such that for any $0< i< n$ and any diagram
\[
\xymatrix{
\Lambda^{n}_i \ar[r] \ar[d] & S \\
\Delta^n \ar@{..>}[ru] & 
}
\]
of solid arrows, there exists a dotted arrow filling the diagram.
Here $\Lambda^{n}_i$ is the $i$-th horn and $\Delta^n$
is the standard $n$-simplex.
Following \cite{HTT}, in the sequel
we call quasi-categories {\it $\infty$-categories}.
A functor of $\infty$-categories $S\to S'$ is
a map of simplicial sets.
By the definition, $\infty$-categories 
form a full subcategory of simplicial sets.
It contains Kan complexes.
The $\infty$-categories also generalize (nerves of) ordinary categories (cf. \cite[1.1.2]{HTT}).
Let $\Delta^1$ be the standard $1$-simplex. It can be regarded as
the nerve of the category $\{0,1\}$ which consists of two objects $0$, $1$,
and the nondegenerate morphism $0\to 1$.
Similarly, $\Delta^0$ can be considered to be the category having one object
with the identity.
Let $S$ be the nerve of the category $\{ 0\stackrel{a,b}{\rightleftarrows} 1\}$
such that $a\circ b=\textup{Id}, b\circ a=\textup{Id}$.
Three simplicial sets $\Delta^1$, $\Delta^0$ and $S$ are all {\it weak homotopy equivalent} to one another,
and they are {\it not isomorphic} to one another as simplicial sets.
However, from the viewpoint of category theory,
we should consider that $\Delta^0$ is ``equivalent'' to $S$,
and $\Delta^1$ is not ``equivalent'' to $\Delta^0$ and $S$.
Hence it is necessary to have a correct notion of
equivalences which generalizes 
the notion of equivalences of ordinary categories.
The important concept we first recall is {\it categorical equivalences} between
simplicial sets.
Let $\textup{Set}_{\Delta}$ be the category of
simplicial sets and let $\Catdel$ be the category of simplicial categories,
in which morphisms are simplicial functors.
Here a simplicial category is a category enriched over the
category of simplicial sets.
Let $\mathcal{H}$ be the homotopy category of ``spaces'',
that is, the category obtained from $\Setdel$ by inverting
weak homotopy equivalences.
To a simplicial category $\mathcal{C}$, applying $\Setdel\to \mathcal{H}$
to the mapping complexes in $\mathcal{C}$ we associate
an $\mathcal{H}$-enriched category $\textup{h}\mathcal{C}$.
Let $\mathcal{C}$ and $\mathcal{D}$ be two simplicial categories.
A simplicial functor $F:\mathcal{C}\to \mathcal{D}$
is said to be an {\it $($Dwyer-Kan$)$ equivalence}
(resp. {\it essentially surjective})
if the induced functor
$\textup{h}\mathcal{C}\to \textup{h}\mathcal{D}$ is an
equivalence of $\mathcal{H}$-enriched categories
(resp. essentially surjective).
A simplicial functor $F:\mathcal{C}\to \mathcal{D}$
is fully faithful if, for any two objects $C,C'\in \mathcal{C}$
the induced morphism $\Map_{\mathcal{C}}(C,C')\to \Map_{\mathcal{D}}(F(C),F(C'))$ is a weak homotopy equivalence.
There is an adjoint pair \cite[1.1.5]{HTT}
\[
\mathfrak{C}:\Setdel\leftrightarrows \Catdel:\textup{N}.
\]
In this paper we will not use the detailed constructions of this adjunction and refer to \cite[1.1]{HTT} for the definition of $\mathfrak{C}$ and $\textup{N}$,
but an important point is that the pair is a Quillen equivalence
with respect to suitable model structures (see below).
The functor $\textup{N}$ is called the simplicial nerve functor.
In fact, if $\mathcal{C}$ is an ordinary category
regarded as a simplicial category, then
$\textup{N}(\mathcal{C})$ coincides with
the usual nerve, and thus
the simplicial nerve functor
generalizes the classical nerve functor
to the $\infty$-categorical setting.
A map of simplicial sets $F:S\to T$ is a {\it categorical equivalence}
(resp. {\it essentially surjective, fully faithful})
if $\mathfrak{C}(F):\mathfrak{C}(S)\to \mathfrak{C}(T)$ is
an equivalence (resp. essentially surjective, fully faithful)
of simplicial categories.

For a simplicial set $S$
we define $\textup{h}(S)$
to be $\textup{h}(\mathfrak{C}(S))$.
Here we ignore the $\mathcal{H}$-enrichment of $\textup{h}(\mathfrak{C}(S))$
and refer to $\textup{h}(S)$ as a homotopy category of $S$.
Here we will describe an alternative construction of
a homotopy category of $S$ when $S$ is an $\infty$-category \cite[1.2.3]{HTT}.
Let $\pi(S)$ be the category defined in the following way.
The objects of $\pi(S)$ are the vertices of $S$.
For $f:\Delta^1\to S$, $f(\{0\})$ and $f(\{1\})$
are said to be the source and the target, respectively.
Let $s,s'\in S$ be two objects and let $f,g:\Delta^1\rightrightarrows S$
be edges.
Suppose that $f$ and $g$ have the same source $s$
and target $s'$.
Then $f$ and $g$ is said to be homotopic if there exists
$\Delta^2\to S$ determined by
\[
\xymatrix{
s \ar[r]^f \ar[rd]_g& s' \ar[d]^{\textup{Id}}\\
  &  s' .
}
\]
Then the relation of homotopy is an equivalence relation on edges
from $s$ to $s'$.
Let $\Hom_{\pi(S)}(s,s')$ be the set of homotopy classes of edges
joining $s$ to $s'$.
Using the definition of $\infty$-categories,
we can define a composition law on the homotopy classes of edges.
This yields a category $\pi(S)$ which turns out to be equivalent to $\textup{h}(\mathfrak{C}(S))$.
Abusing notation we often write $\textup{h}(S)$ for $\pi(S)$.

The pair of adjoint $(\mathfrak{C},\textup{N})$ plays an important role in the various constructions
of $\infty$-categories and their functors and so on.
For various applications,
it is better to view this adjoint
as a Quillen adjoint pair with respect to suitable
model structures on $\Setdel$ and $\Catdel$ rather than a usual
adjoint pair. 
The category $\Setdel$ admits a model structure,
in which a weak
equivalence is a categorical equivalence, and a cofibration is
a monomorphism (\cite{J}, \cite[2.2.5.1]{HTT}).
Our reference of model categories are \cite{Hov98} and \cite[Appendix]{HTT}.
It turns out that an object is fibrant in this model category
if and only if
it is an $\infty$-category. We refer to this model structure
as {\it Joyal model structure}.
There exists a model structure on $\Catdel$
such that the weak equivalences are equivalences
and fibrant objects are simplicial categories
whose mapping complexes are Kan complexes
(see for the details \cite{Berg}, \cite[A 3.2]{HTT}).
Then the adjunction $(\mathfrak{C},\textup{N})$
is a Quillen equivalence with respect to these model structures
(see \cite[2.2.5.1]{HTT}).
For example, we can use this adjoint as follows.
Let $\mathsf{M}$ be a simplicial model category.
Then the full subcategory $\mathsf{M}^{\circ}$ spanned by
cofibrant-fibrant objects is a fibrant simplicial category.
(For a model category $\mathsf{M}$ we shall denote by
$\mathsf{M}^{\circ}$ the full subcategory spanned by cofibrant-fibrant objects.)
Applying the simplicial nerve functor to
$\mathsf{M}^{\circ}$ we obtain an $\infty$-category $\textup{N}(\mathsf{M}^{\circ})$. 
There is another method by which one can obtain the $\infty$-category
from a model category.
Let $\mathsf{M}$ be a model category
and let $\mathsf{M}^{c}$ be a full subcategory spanned by cofibrant
objects. Let $W$ be the collection of weak equivalences
in $\mathsf{M}^{c}$.
Then there exist an $\infty$-category
$\NNNN(\mathsf{M}^{c})[W^{-1}]$ and a functor $\NNNN(\mathsf{M}^{c})\to\NNNN(\mathsf{M}^{c})[W^{-1}]$ such that for any $\infty$-category $\mathcal{C}$
the composition induces a fully faithful functor
\[
\Fun(\NNNN(\mathsf{M}^{c})[W^{-1}],\mathcal{C})\to \Fun(\NNNN(\mathsf{M}^{c}),\mathcal{C})
\]
and its essential image is spanned by functors
$\NNNN(\mathsf{M}^{c}) \to \mathcal{C}$ which carries endges in $W$ to equivalences in $\mathcal{C}$ (see \cite[1.3.4]{HA}).
If $\mathsf{M}$ is a simplicial model category,
then there exists a natural categorical equivalence
$\NNNN(\mathsf{M}^\circ) \simeq \NNNN(\mathsf{M}^{c})[W^{-1}]$.
We call
$\NNNN(\mathsf{M}^{c})[W^{-1}]$
the underlying $\infty$-category of $\mathsf{M}$.
When $\mathsf{M}$ is a combinatorial model category, then
$\NNNN(\mathsf{M}^{c})[W^{-1}]$ is presentable (see \cite[5.5]{HTT}
for the notion of presentable $\infty$-categories).
If $\mathsf{M}$ is a symmetric monoidal model category,
$\NNNN(\mathsf{M}^{c})[W^{-1}]$ inherits a symmetric monoidal
structure (see \cite[4.1.3]{HA}), which we shall refer to as the underlying
symmetric monoidal $\infty$-category.

The Joyal model structure on $\Setdel$ is relevant to the usual
model structure introduced by Quillen (\cite{Q}), in which a weak equivalence
is a {\it weak homotopy equivalence}.
According to \cite{J} \cite[2.2.5.8]{HTT}, we have the following implications
\[
(\textup{isomorphisms})\Rightarrow(\textup{categorical equivalences})\Rightarrow (\textup{weak homotopy equivaleces}).
\]

Let $S$ be a simplicial set.
An {\it object} in $S$ is a vertex $\Delta^0\to S$.
A {\it morphism} in $S$ is an edge $\Delta^1\to S$, and
when $S$ is an $\infty$-category, a morphism $\Delta^1\to S$
is said to be an {\it equivalence} if it gives rise to an isomorphism
in the homotopy category $\textup{h}S$.
Let $\mathcal{C}$ and $\mathcal{D}$ be two $\infty$-categories.
Define $\Fun(\mathcal{C},\mathcal{D})$ to be
the simplicial set $\Map_{\textup{Set}_{\Delta}}(\mathcal{C},\mathcal{D})$
which parametrizes maps from $\mathcal{C}$ to $\mathcal{D}$,
that is, a map $\Delta^n\to \Fun(\mathcal{C},\mathcal{D})$
amounts to a map $\mathcal{C}\times\Delta^n\to \mathcal{D}$.
By \cite[1.2.7.3]{HTT}, the simplicial set $\Fun(\mathcal{C},\mathcal{D})$
is an $\infty$-category.
We shall refer to an object of $\Fun(\mathcal{C},\mathcal{D})$
as a {\it functor} from $\mathcal{C}$ to $\mathcal{D}$.
We shall refer to a {\it morphism} (resp. an {\it equivalence}, i.e., a {\it morphism} which induces an isomorphism in $\textup{h}\Fun(\mathcal{C},\mathcal{D})$) in $\Fun(\mathcal{C},\mathcal{D})$
as a {\it natural transformation} (resp. a {\it natural equivalence}).
We define $\Map(\mathcal{C},\mathcal{D})$ to be the largest
Kan complex of $\Fun(\mathcal{C},\mathcal{D})$.
Namely, $\Map(\mathcal{C},\mathcal{D})$ is the subcategory spanned by
natural equivalences.
Define $\Catinf^{\Delta}$ to be a fibrant simplicial category
whose objects are small $\infty$-categories, and
whose hom simplicial set (between $\mathcal{C}$ and $\mathcal{D}$)
is $\Map(\mathcal{C},\mathcal{D})$.
Let $\Catinf$ be the simplicial nerve of $\Catinf^{\Delta}$ (cf. \cite[Chapter 3]{HTT}).
We shall refer to $\Catinf$ as the $\infty$-category of (small) $\infty$-categories.
We shall denote by $\widehat{\mathcal{C}\textup{at}}_{\infty}$ the $\infty$-category of
(large) $\infty$-categories.
We let $\mathcal{S}$ be the full subcategory of
$\Catinf$ spanned by Kan complexes, which we shall refer to
as the $\infty$-category of spaces.

Let $S$ be an $\infty$-category and let $s,s'$ be two objects in $S$.
In the course of the paper
 we sometimes discuss
``the mapping space'' from $s$ to $s'$.
The direct way to the definition is
to define $\Map_S(s,s')$ to be the complex $\Map_{\mathfrak{C}(S)}(s,s')$.
Remembering the relationship between $\infty$-categories and simplicial
categories, we should regard the simplicial set $\Map_S(s,s')$
as an object in the homotopy category $\mathcal{H}$ of spaces.
The simplicial set $\Map_{\mathfrak{C}(S)}(s,s')$ equips associative
compositions (varying $s$ and $s'$), but it is not a
Kan complex in general.
There are several ways to construct a simplicial set that represents
the weak homotopy type of $\Map_{\mathfrak{C}(S)}(s,s')$.
For example, the Kan complex of left morphisms
$\Hom^{\textup{L}}_S(s,s')$
determined by
\[
\Hom_{\Setdel}(\Delta^n,\Hom^{\textup{L}}_S(s,s'))=\{f:\Delta^{n+1}\to S\ |\ 
f|_{\Delta^{\{0\}}}=s,\ f|_{\Delta^{\{1,\ldots,n+1\}}}\ \textup{is constant at }s'\}
\]
represents the weak homotopy type of
 $\Map_{\mathfrak{C}(S)}(s,s')$.
(see \cite[1.2.2]{HTT} for
details of mapping complexes).

\subsection{$\infty$-category of quasi-coherent complexes}\label{subsec:2-2}
We refer to \cite{LM} as the general reference of the notion
of algebraic stacks.
In this paper, all algebraic stacks (and schemes)
are assumed to be {\it Deligne-Mumford}, {\it quasi-compact} and to {\it have affine diagonal}.
We fix three Grothendieck universes $\mathbb{U}_1\in \mathbb{U}_2\in \mathbb{U}_3$ such that $\mathbb{U}_1$ contains all finite ordinals.
We suppose that all schemes, rings and others belong to $\mathbb{U}_1$
and all (pre)sheaves are $\mathbb{U}_1$-small.
Entries in $\mathbb{U}_1$ (resp. $\mathbb{U}_2$, $\mathbb{U}_3$)
are called small (resp. large, super-large).
By a vector bundle on an algebraic stack $\XX$
we mean a locally free $\OO_{\XX}$-module of finite type.

In the rest of this section, $R$ is a commutative ring.
We denote by $\Coch(R)$ the category of cochain complexes of $R$-modules.
Let $\chain(R)$ be the category of chain complexes of $R$-modules.

Let $S=\Spec A$ be an affine scheme.
A quasi-coherent $\OO_S$-complex can be regarded as a set of data $\{M_B, \alpha\}_{\Spec B\to \Spec A}$ consisting
of an (unbounded) complex $M_B$ of $B$-modules
for any $\Spec B\to \Spec A$,
and an isomorphism
$\alpha_{\phi}:M_{B}\otimes_BB'\to M_{B'}$ for any $\phi:\Spec B'\to \Spec B$ over $\Spec A$, such that $\alpha_\phi$'s satisfy the cocycle condition, i.e.,
$\alpha_{id}={id}$ and
for each $\phi\circ \psi:\Spec B''\to \Spec B'\to\Spec B$ over $\Spec A$
we have $\alpha_{\phi\circ \psi}=\alpha_{\psi}\circ\alpha_{\phi}$.
Let $\textup{QC}(S)$ be the category of quasi-coherent $\OO_S$-complexes.
By the projective model structure \cite[Section 4.2]{Hov98},
$\textup{QC}(S)$ forms a symmetric monoidal model category, in which
the weak equivalences are quasi-isomorphisms, and fibrations are
termwise surjections.
Let $\textup{QC}(S)^{c}\subset \textup{QC}(S)$ be the full subcategory spanned by cofibrant objects.
We define the stable $\infty$-category $\DQ(S)$ of quasi-coherent complexes
to be the underlying $\infty$-category (see Section \ref{subsec:2-1}).
For the generalities of stable $\infty$-categories,
see \cite[Chapter 1]{HA}.
Let $J$ be the category of affine schemes over $\XX$
and we abuse notation and often write $J$ for the nerve $\textup{N}(J)$.
A marked simplicial set is a pair $(X,E)$ consisting
of a simplicial set $X$ and a set $E$ of edges of $X$ that includes all
degenerate edges. A morphism $(X,E)\to (X',E')$ of marked simplicial
sets is
a map $f:X\to X'$ such that $f(E)\subset E'$.
Let
$\widehat{\textup{Set}}_{\Delta}^{+}$ denote the
category of marked large simplicial sets \cite[3.1]{HTT},
which is endowed with a simplicial combinatorial
model structure \cite[3.1.3.7]{HTT}.
Let $(\widehat{\textup{Set}}_{\Delta}^{+})^{\circ}$ denote the full subcategory
spanned by cofibrant-fibrant objects.
There exists a natural categorical equivalence $\textup{N}((\widehat{\textup{Set}}_{\Delta}^{+})^{\circ})\simeq \widehat{\mathcal{C}\textup{at}}_{\infty}$.
Let $J^{op}\to \widehat{\textup{Set}}_{\Delta}^{+}$ be a functor
which sends $S\to \XX\in J$ to $(\NNNN(\textup{QC}(S)^{c}),W_S)$ and sends $f:S'\to S$ (over $\XX$)
to $f^*:(\NNNN(\textup{QC}(S)^{c}),W_S)\to (\NNNN(\textup{QC}(S')^{c}),W_{S'})$.
Here $W_S$ and $W_{S'}$ are collections of weak equivalences in
$\NNNN(\textup{QC}(S)^{c})$ and $\NNNN(\textup{QC}(S')^{c})$ respectively.
We can associate the underlying $\infty$-category
$\NNNN(\textup{QC}(S)^{c})[W_{S}^{-1}]$ to $(\NNNN(\textup{QC}(S)^{c}),W_S)$
in a functorial fashion by using the functorial
fibrant replacement in $\widehat{\textup{Set}}_{\Delta}^{+}$.
Namely, we have $J^{op}\to \widehat{\textup{Set}}_{\Delta}^{+}$
which carries $S\to\XX$ to $\NNNN(\textup{QC}(S)^{c})[W_{S}^{-1}]$.
It gives rise to a functor
$J^{op}\to \widehat{\mathcal{C}\textup{at}}_{\infty}$.
Following \cite{TV} and \cite{BFN}, for an algebraic stack $\XX$
we define the stable $\infty$-category
$\DQ(\XX)$ of quasi-coherent complexes
by
\[
\DQ(\XX):=\lim_{S\to \XX \in J}\DQ(S)
\]
where $\lim$ means a limit in the $\infty$-category $\widehat{\mathcal{C}\textup{at}}_{\infty}$.
For a morphism $f:S\to \XX$,
we define $f^*:\DQ(\XX)\to \DQ(S)$
to be the natural projection $\lim_{S\to \XX}\DQ(S)\to \DQ(S)$.
Since $\DQ(S)$ is a presentable $\infty$-category for any affine scheme
$S$, a standard cardinality estimation shows
that $\DQ(\XX)$ is presentable when the stack (as a functor)
$\mathcal{X}:\textup{CAlg} \to \mathcal{S}$
is an accessible functor, that is,
there exists a homotopy left cofinal small full subcategory of $J^{op}$
(consider $\kappa$-compact objects in $\textup{CAlg}$
for some cardinal $\kappa$).
(This condition holds for our schemes and algebraic stacks in Section \ref{subsec:2-3}.)
Here $\textup{CAlg}$ is the $1$-category of (usual) commutative rings.

Let $\lcatinf^{\textup{sMon}}$ be the $\infty$-category of
symmetric monoidal $\infty$-categories whose morphisms
are symmetric monoidal functors.
Since the notion of symmetric monoidal $\infty$-categories
will not appear until Section \ref{sec:5}, we postpone giving the definitions of
these notions.
We quickly define the symmetric monoidal $\infty$-category
$\DQ^\otimes(\XX)$ (without giving precise definitions concerning
symmetric monoidal $\infty$-categories).
Since $\textup{QC}(S)$ are symmetric monoidal model categories,
by \cite[4.1.3.4]{HA} the above
functor $J^{op}\to \lcatinf$ is promoted to a functor
$J^{op}\to \lcatinf^{\textup{sMon}}$.
The symmetric monoidal $\infty$-category $\DQ^\otimes(\XX)$ is defined to
be a limit of $J^{op}\to \lcatinf^{\textup{sMon}}$.
The underlying $\infty$-category is equivalent to
$\DQ(\XX)$ defined above.
The symmetric monoidal structure on $\DQ(\XX)$
induces a symmetric monoidal structure on the homotopy category
$\textup{D}_{\textup{qcoh}}(\XX)=\textup{h}\DQ(\XX)$.

When $f:\XX\to \Spec R$ is an algebraic stack over $R$,
$\DQ(\XX)$ has an $R$-linear structure.
(In this paper, the notion of $R$-linear $\infty$-categories is not needed except
the application of derived Morita theory. Hence the reader who is not familiar
with $\infty$-operads may skip this paragraph.)
There exists the pullback functor
$f^*:\DQ^{\otimes}(R)\to \DQ^\otimes(\XX)$.
It induces an ``action'' of
$\DQ^{\otimes}(R)$ on $\DQ(\XX)$.
Roughly speaking, this action
is determined by the composite
\[
\DQ(R)\times \DQ(\XX)\stackrel{f^*\times \textup{id}}{\longrightarrow} \DQ(\XX)\times \DQ(\XX)\stackrel{\otimes}{\longrightarrow} \DQ(\XX)
\]
and homotopy coherence of associativities.
Let $\mathcal{M}\textup{od}_{\DQ(R)}$ be
the $\infty$-category of left module objects in $\widehat{\mathcal{C}\textup{at}}_{\infty}$ over the monoidal $\infty$-category $\DQ(R)$
(see \cite[3.3.3, 4.2]{HA}).
We refer to an object as an $R$-linear $\infty$-category.
The above action exhibits $\DQ(\XX)$ as an $R$-linear $\infty$-category.
By a straightforward construction from $J^{op}\to \lcatinf^{\textup{sMon}}$
via the machinery of $\infty$-operads \cite{HA},
we have a natural functor
\[
J^{op}\longrightarrow \mathcal{M}\textup{od}_{\DQ(R)}
\]
which extends $J^{op}\to \widehat{\mathcal{C}\textup{at}}_{\infty}$.
Let $\mathcal{P}r^{L}$ be the subcategory
of $\widehat{\mathcal{C}\textup{at}}_{\infty}$ spanned by
presentable $\infty$-categories, in which functors are left adjoints
(see \cite[5.5.3]{HTT}). 
The $\infty$-category $\mathcal{P}r^{L}$ inherits the (symmetric)
monoidal structure
described in \cite[6.3.1.14]{HA}.
Let $\mathcal{M}\textup{od}_{\DQ(R)}(\mathcal{P}r^L)$ denote
the $\infty$-category of left module objects in $\mathcal{P}r^L$
over the monoidal $\infty$-category $\DQ(R)$.
Since the left module map $\DQ(R)\times \DQ(S)\to \DQ(S)$
is a map which preserves colimits separately in each variable,
thus $J^{op}\to \mathcal{M}\textup{od}_{\DQ(R)}$
yields $J^{op}\to \mathcal{M}\textup{od}_{\DQ(R)}(\mathcal{P}r^L)$.
Take a limit of $J^{op}\to \mathcal{M}\textup{od}_{\DQ(R)}(\mathcal{P}r^L)$. According to \cite[4.2.3.3]{HA},
the underlying category of this limit is equivalent to $\DQ(\XX)$.

\vspace{1mm}

\subsection{Schemes and stacks}\label{subsec:2-3}
Let $k=R$ be a field.
Let $\XX$ be an algebraic stack over $k$.
In the main results of
this paper we will treat the following two cases:
\begin{enumerate}
\renewcommand{\labelenumi}{(\roman{enumi})}

\item $\XX$ is a
noetherian scheme which has a very ample invertible sheaf
(e.g., quasi-projective varieties).

\item $\XX$ is a tame separated (Deligne-Mumford) algebraic stack of the form $[X/G]$
where $X$ is a finitely generated
noetherian scheme and $G$ is a linear algebraic group
acting on $X$. Suppose further that the
coarse moduli space is quasi-projective and $X$ has
a $G$-ample invertible sheaf.

\end{enumerate}

\begin{Remark}
The following are examples of algebraic stacks which satisfy the condition (ii).
\begin{enumerate}

\item GIT stable quotients whose stabilizer groups are all finite group. Let $X$ be a separated scheme of finite type over a field, endowed with action of a linearly reductive
group $G$. Assume that there exists a $G$-linearized ample invertible sheaf $\mathcal{L}$
and $X^s(\mathcal{L})$ is the open subset of stable points. Suppose furthermore that every stabilizer is finite.
Then the quotient stack $[X^s(\mathcal{L})/G]$ satisfies the condition (ii).
Such quotients often arise from Geometric Invariant Theory.

\item Separated and smooth tame Deligne-Mumford stacks
which satisfy the conditions:
(a) it is generically a scheme,
(b) the coarse moduli space is quasi-projective (see \cite[Theorem 1.2]{Tot}).
(For example, toric stacks (orbifolds) cf. \cite{Iwa}.)

\item Moduli stacks of stable curves, stable maps, polarized abelian varieties, Calabi-Yau manifolds in characteristic zero.

\end{enumerate}

\end{Remark}

Let us recall the notion of {\it perfect stacks}
introduced in \cite{BFN}.
In loc. cit., the authors offer us
the concept in
the framework of derived stacks
and prove derived Morita theory for perfect stacks,
but here we consider only usual
algebraic stacks.
Let $\XX$ be an algebraic stack.
Let $\DP(\XX)\subset \DQ(\XX)$ be the full subcategory
consisting of perfect complexes.
(A strictly perfect complex on $\XX$ is
a bounded complex of vector bundles.
A complex is said to be perfect if locally on the smooth site of $\XX$
it is quasi-isomorphic to a strictly perfect complex.
According to \cite[3.6, 3.7]{Kri} under the assumption
of (i) and (ii) any perfect complex
of $\OO_{\XX}$-modules is quasi-isomorphic to
a perfect complex of quasi-coherent sheaves.)
An algebraic stack $\XX$ is said to be perfect
if the $\infty$-category $\textup{Ind}\DP(\XX)$
of Ind-objects \cite[5.3]{HTT} of perfect complexes is naturally equivalent to
$\DQ(\XX)$.
A large class of stacks satisfies perfectness (e.g. quasi-compact
and separated schemes, quotient stacks in characteristic zero, algebraic stacks satisfying (i) or (ii), etc, see \cite{BFN}, \cite{Azu},
Corollary~\ref{pffincd}). 
If $\XX$ is a perfect stack, then
$\DP(\XX)$ is the full subcategory spanned by compact objects in $\DQ(\XX)$ on one hand,
and $\DQ(\XX)$ is $\textup{Ind}\DP(\XX)$ on the other hand.
Consequently, we can transform $\DP(\XX)$
into $\DQ(\XX)$ and transform $\DQ(\XX)$ into $\DP(\XX)$
in the categorical fashion.
In particular, we can consider an exact functor $\DP(\XX) \rightarrow \DP(S)$ 
to be a colimit-preserving functor $\DQ(\XX) \rightarrow \DQ(S)$, which preserves
compact objects.

\begin{Lemma}
Let $\XX$ be a tame Deligne-Mumford stack which is separated and of finite
type over $\ZZ$.
Suppose that its coarse moduli space is a scheme.
Then compact and dualizable objects in $\dqc(\stk X)$ coincide
$($see Section \ref{sec:4} for the notion of dualizable objects$)$.
\end{Lemma}

\Proof
To see that dualizable objects are compact,
it is enough to show that the (derived) global section functor $\Gamma(\stk X,-)$ preserves colimits
since the functor $\Hom(P,-)$ is equivalent to $\Gamma(\stk X,P^*\otimes -)$
for any dualizable object $P$ and the functor
 $P^*\otimes -$ preserves colimits. Here $P^*$ is a dual of $P$.
By our assumption on $\stk X$, we have
a coarse moduli space $p:\stk X \rightarrow M$
such that $M$ is quasi-compact and separated.
Thus $M$ is a perfect stack (cf. \cite[Section 3]{BFN}).
Notice that
the dualizable object $\OO_M$ is compact in $\DQ(M)$.
Since $\OO_M$ is compact, the functor $\Gamma(M,-)$ preserves colimits.
Hence to see that $\Gamma(\stk X,-)$ preserves colimits,
it is sufficient to show that the pushforward $p_*$ preserves colimits.
There exist an \'etale surjective morphism $U\to M$
and a Cartesian diagram
\begin{eqnarray*}
\begin{xymatrix}
{
[W/G'] \ar[r]\ar[d]^{p_U} & \stk X \ar[d]^p\\
U \ar[r] & M\\
}
\end{xymatrix}
\end{eqnarray*}
where $U$ is an affine scheme
and $[W/G']$ is a quotient stack of a finite scheme $W$ (over $U$)
by action of a finite (\'etale) group scheme $G'$ over $U$. 
Then $[W/G']$ is perfect by \cite[Proposition 3.26]{BFN}
(we here use the tameness of $\XX$).
Thus by \cite[Proposition 3.10, 3.23]{BFN} $p_U$ is a perfect morphism and
$p_{U*}$ preserves small colimits.
Since $U \rightarrow M$ is \'etale surjective, (using descent and base change theorem)
we see that $p_*$ 
also preserves small colimits.
Conversely, to see that compact objects are dualizable,
it is enough to repeat the same argument in the proof of \cite[Lemma 3.20]{BFN}
for $p:\stk X \rightarrow M$ and the affine covering map $U\rightarrow M$.
\QED

According to \cite[Proposition 3.9]{BFN} an algebraic stack $\XX$ is
perfect if and only if
compact and dualizable objects in $\dqc(\stk X)$ coincide
and $\dqc(\stk X)$ is compactly generated.
The recent powerful result of the existence of compact generators
by To\"en show that a separated and quasi-compact Deligne-Mumford
stack has a compact generator if its coarse moduli space is a scheme
(see \cite[4.2]{Azu}).
Thus we have:

\begin{Corollary}\label{pffincd}
Let $\XX$ be a tame Deligne-Mumford stack which is separated and of finite
type over $\ZZ$.
Suppose that a coarse moduli space for $\XX$ is a scheme.
Then $\stk X$ is a perfect stack.
\end{Corollary}

\section{Extension Lemmas}\label{sec:3}

In this section let $\XX$ and $S$ be perfect stacks.
Let $\dvb(\stk X)$ (resp. $\dvb(S)$) denote the full subcategory of $\dpf(\stk X)$ (resp. $\dpf(S)$),
spanned by quaisi-coherent complexes which are quasi-isomorphic to
vector bundles placed in degree zero
on $\stk X$ (resp. $S$).

\begin{Lemma}
\label{ff;1cat}
We have the followings:
\begin{enumerate}
\renewcommand{\labelenumi}{(\roman{enumi})}
\item $\DV(S)$ is equivalent to a $1$-category $($cf. \cite[Section 2.3.4]{HTT}$)$.

\item Let $\mathcal{E}$ be an $\infty$-category. The functor $\Fun (\textup{N}(\textup{h}\mathcal{E}),\DV(S))\to \Fun(\mathcal{E},\DV(S))$
associated to the projection $\mathcal{E}\to \textup{N}(\textup{h}\mathcal{E})$
is a categorical equivalence.

\end{enumerate}
\end{Lemma}

\Proof
We first prove (i).
Note that
for any locally free sheaves $E$
and $F$ on $S$, the Ext-group $\Ext^i(E,F)$
is zero for $i<0$.
It follows that for
every pair of objects $E,F\in \DV(S)$,
the mapping space $\Map_{\DV(S)}(E,F)$
is discrete, that is, 0-truncated. Therefore $\DV(S)$ is equivalent to
a 1-category
(cf. \cite[2.3.4.18]{HTT}).
The claim (ii) follows from \cite[2.3.4.12]{HTT}.
\QED

The homotopy category $\textup{h}\DV(\XX)$ is a 1-category
whose objects are vector bundles
on $\XX$, placed on degree zero.
A morphism $E\to F$
in $\textup{h}\DV(\XX)$ can be considered to be a morphism
of locally free sheaves on $\XX$.

\vspace{2mm}

\begin{Lemma}
\label{bundlecolimit}
For any $n\geq 0$,
let $I_n$ denote the simplicial set defined as follows:
\begin{eqnarray*}
I_0 &:=& \Delta^0\\
I_{n+1} &:=&
(I_n\times\Delta^1) \coprod_{I_n\times \Delta^{\{0\}}}
(I_n \times\Delta^1).
\end{eqnarray*}
Let $P\in\DQ(\XX)$
be a strict perfect complex on $\XX$ which lies in $(-\infty,0]$.
Suppose that $P$ is represented by the complex of the form
\[
\ldots 0 \to P^{-n}\to \ldots \to P^{-1}\to P^0\to 0\to \ldots
\]
where $P^{i}$ is a vector bundle placed in degree $i$.
Then there exists a set of diagrams 
$\{p_k:I_k \rightarrow \DQ(\XX)\}_{k\geq 0}$
such that
the complex $(\sigma^{\le -n+k}P)[-n+k]$ is a colimit of $p_k$
for any $k\geq 0$
and the restriction $p_k|_{I_{k-1}\times\Delta^0}$ is $p_{k-1}$
for any $k\geq 1$.

\end{Lemma}

\Proof
We will inductively construct $p_k$.
Let $p_0$ be the map $\Delta^0\rightarrow \DQ(\XX)$
which sends $0\in\Delta^0$ to $P^{-n}$.
Now suppose that we have constructed $p_k$ for any $k\leq l$.
Let $P'$ denote the complex $(\sigma^{\le -n+l}P)[-n+l]$.
Since $P'$ is a colimit of $p_l$,
the canonical morphism of complexes
$P' \rightarrow P^{-n+l+1}$
induces a map $q: I_k\times\Delta^1 \to \DQ(\XX)$
such that $q|_{I_k\times\Delta^{\{0\}}}=p_l$
and $q|_{I_k\times\Delta^{\{1\}}}$ is a constant diagram
with value $P^{-n+l+1}$.
Let $q': I_k\times\Delta^1 \to \DQ(\XX)$ be a map
such that $q'|_{I_k\times\Delta^{\{0\}}}=p_l$
and $q'|_{I_k\times\Delta^{\{1\}}}$ is a constant diagram
with value $0$.
Then we obtain a diagram $p_{l+1}:I_{l+1}\rightarrow \DQ(\XX)$
by gluing $q$ and $q'$ along $p_l$.
The (homotopy) pushout $0 \leftarrow P'
\rightarrow P^{-n+l+1}$ is a colimit of $p_{l+1}$
by \cite[4.4.2.2]{HTT}.
Consider the mapping cylinder 
\[
\xymatrix@R=7mm @C=5mm{
 \ldots \ar[r] & 0 \ar[r] \ar[d]& P^{-n} \ar[r]^{(-1)^{n-l}d_P}\ar[d]& \ldots  \ar[r]& P^{-n+l-1} \ar[r]^{(-1)^{n-l}d_P}\ar[d]& P^{-n+l}\ar[d] \ar[r] & 0\\
\ldots \ar[r] & P^{-n} \ar[r]& P^{-n}\oplus P^{-n+1} \ar[r] & \ldots\ar[r] & P^{-n+l-1}\oplus P^{-n+l} \ar[r] & P^{-n+l}\oplus P^{-n+l+1} \ar[r] &0 
}
\]
of $P'\to P^{-n+l+1}$.
Here we regard $P^{-n+l+1}$ as a complex whose degree zero term is $P^{-n+l+1}$.
Let $P''$ denote the lower complex.
The vertical arrows in the mapping cylinder
are split monomorphisms and thus
$(\sigma^{\le -n+l+1}P)[-n+l+1]$ is a homotopy colimit of the diagram $0\leftarrow P'\to P''$
\QED

\begin{Remark}
\label{sketo}
An analogues result holds for a 
bounded complex of quasi-coherent sheaves $P^{\bullet}$
such that $P^i=0$ for $i>0$.
\end{Remark}

\begin{Lemma}
\label{inductive}
Let $\XX$ be an algebraic stack which satisfies either the condition (i)
or (ii) in Section \ref{subsec:2-3}.
Let $P$ be a complex of quasi-coherent sheaves, i.e., $P\in\DQ(\XX)$. Then
there exist a filtered system of complexes $\{E(n,m)\}_{n\ge0,m\ge0}$
and a quasi-isomorphism
$\lim_{n,m} E(n,m)\to P$ such that
$E(n,m)$ is quasi-isomorphic to $\sigma^{\ge -m}\tau^{\le n}P$,
$E(n,m)$ is a complex which in each degree is an infinite direct sum
of invertible sheaves, and $E(n,m)^i$ is zero for $i>n$ and $i<-m$.
\end{Lemma}

\Proof
It follows from \cite[2.3.2]{TT} and its proof.
\QED

Let $\DP^{\le 0}(\XX)$
denote the full subcategory of $\DP(\XX)$
spanned by complexes 
which are equivalent to complexes $P^{\bullet}$
such that if $i\le 0$, $P^{i}$ is a vector bundle,
and if $i>0$ then $P^i=0$.
The theory of left Kan extensions \cite[4.3.2.16]{HTT} for $\infty$-categories
provides the following lemma:

\begin{Lemma}
\label{leftkan}
Let $\kappa:\Fun(\DP^{\le0}(\XX),\DQ(S))\to \Fun(\DV(\XX),\DQ(S))$
be the functor induced by the inclusion $\DV(\XX)\subset \DP^{\le0}(\XX)$.
Let $\mathcal{K}'\subset \Fun(\DV(\XX),\DQ(S))$ be the full
subcategory spanned by functors $\DV(\XX)\to \DQ(S)$
whose essential images lie in $\DV(S)$.
Let $\mathcal{K}\subset \Fun(\DP^{\le0}(\XX),\DQ(S))$
be the full subcategory spanned by the functors
$\Phi:\DP^{\le0}(\XX)\to \DQ(S)$ such that $\Phi$ is a
left Kan extension of $\Phi|_{\DV(\XX)}:\DV(\XX)\to \DQ(S)$ and
$\kappa(\Phi)\in \mathcal{K}'$.
Then the induced map $\mathcal{K}\to \mathcal{K}'$
is a categorical equivalence.
\end{Lemma}

$\bullet$ In what follows we will assume that
$\XX$ has the {\it resolution property}, that is, every coherent sheaf $\mathcal{F}$
on $\XX$
admits a surjective morphism $\mathcal{E}\to \mathcal{F}$ from a
vector bundle $\mathcal{E}$.
Under the condition (i) and (ii) in Section \ref{subsec:2-3}
(the existence of a $G$-ample invertible sheaf), $\XX$ has 
the resolution property.
However, note that
the resolution property is not needed in
Lemma~\ref{extensionstep1}, 
\ref{openreduction}.

\vspace{3mm}
We will say that an algebraic stack $\XX$ has cohomological dimension zero
if $\textup{H}^i(\XX,E)$ is zero
for any quasi-coherent sheaf $E$ and $i>0$.

\begin{Proposition}
\label{cofinality}
Suppose that $\XX$ has cohomological dimension zero.
Let $\overline{\Phi}$ be a colimit-preserving functor
$\DQ(\XX)\to \DQ(S)$
such that $\overline{\Phi}(\DP(\XX))$ lies in $\DP(S)$.
Let $\Phi$ be the restriction $\overline{\Phi}|_{\DP^{\le0}(\XX)}:\DP^{\le0}(\XX)\to \DP(S)$
and suppose that $\Phi(\DV(\XX))$ lies in $\DV(S)$.
Then $\Phi$ belongs to $\mathcal{K}$.
\end{Proposition}

\Proof
It is clear that $\kappa(\Phi)$ belongs to $\mathcal{K}'$.
Thus it suffices to prove that $\Phi$ is a left Kan extension
of $\Phi_0:=\Phi|_{\DV(\XX)}:\DV(\XX)\to \DP(S)$.
Recall that $\Phi$ is said to be a left Kan extension
if for any $P\in \DP^{\le0}(\XX)$
the induced functor $p$ in the commutative diagram
\[
\xymatrix{
\DV(\XX)_{/P}\ar[r]\ar[d] & \DP^{\le0}(\XX)\ar[r]^{\Phi} & \DP(S) \\
(\DV(\XX)_{/P})^{\triangleright}\ar[rru]_p & & 
}
\]
is a colimit diagram. (Here the cone point of $(\DV(\XX)_{/P})^{\triangleright}$ maps to $\Phi(P)$.)
To prove this, we may replace $\DP(S)$ by $\DQ(S)$.
Fix a perfect complex $P$ in $\DP^{\le0}(\XX)$. It is quasi-isomorphic to
a strict perfect complex since we impose the resolution property.
Since $\overline{\Phi}$ preserves small colimits,
it is enough to show that $P$ is a colimit of $\DV(\XX)_{/P}\to \DQ(\XX)$.
To this end, let $K=\DV(\XX)$ and take a colimit $R$ of the diagram $K_{/P}\to K\to \DQ(\XX)$.
By Lemma~\ref{bundlecolimit},
$P$ is a colimit of the diagram $p_n:I_n\to \DQ(\XX)$
of vector bundles
(we here use the notation in Lemma~\ref{bundlecolimit}).
Invoking the universality of $P$ and $R$,
we obtain morphisms $P\to R$ and $R\to P$.
Note that the composite
$P\to R\to P$ is equivalent to the identity morphism.
Since $\DQ(\XX)$ is idempotent complete, $R$ has the form $P\oplus P'$,
and $P$ is identified with the direct summand $P\oplus \{0\}\subset R$.
We may and will identify $R$ with $P\oplus P'$.
To complete the proof, it will suffice to prove that $P'$ is
a zero object. Now suppose that $P'$ is not
a zero object. Then there exists $(\theta:E\to P)\in K_{/P}$
such that the corresponding
morphism $\xi:E\to  R$, in the colimit diagram $(K_{/P})^{\triangleright}\to \DQ(\XX)$, whose cone point maps to $R$, induces a non-null-homotopic
morphism $E\to R\simeq P\oplus P'\stackrel{\textup{pr}_2}{\to} P'$.
On the other hand, since cohomological dimension of $\XX$ is zero,
there exists some $E\to P^0$ which represents
$\theta:E\to P$.
Note that the composite $P^0\to P\oplus\{0\}\subset P\oplus P'\stackrel{\textup{pr}_2}{\to}P'$
is null-homotopic.
It follows that $\xi:E\to P\oplus P'$ factors through $E\to P^0\to P$.
Consequently, $E\stackrel{\xi}{\to} P\oplus P'\to P'$ is null-homotopic. It gives rise to
a contradiction, as desired.
\QED

For the ease of notation, in the proofs,
we usually denote by $\CCC$ and $\DDD$ stable presentable $\infty$-categories $\DQ(\XX)$ and
$\DQ(S)$ respectively.
Similarly, we denote by
$\CCC_{\circ}\subset\CCC$ $\DDD_{\circ}\subset \DDD$
the full subcategories consisting of perfect complexes.
(Note that $\textup{Ind}\CCC_{\circ}\simeq\CCC$
and $\textup{Ind}\DDD_{\circ}\simeq\DDD$.)
Let $\CCC_v$ and $\DDD_v$ be the full subcategories of
$\CCC$ and $\DDD$ respectively, spanned by vector bundles (i.e., complexes
which are equivalent to vector bundles).

\begin{Lemma}
\label{extensionstep1}
Let
$\Fun'(\DQ(\XX)^{\times n},\DQ(S))\subset \Fun(\DQ(\XX)^{\times n},\DQ(S))$
be the full subcategory
spanned by functors which preserves colimits in each variable. $($the product $\DQ(\XX)^{\times n}$ is $n$-times $($homotopy$)$ product.$)$
Let $\Fun''(\DP(\XX)^{\times n},\DP(S))\subset \Fun(\DP(\XX)^{\times n},\DP(S))$
be the full subcategory
spanned by functors which preserve finite colimits in each variable.
Then the restriction functor
\[
\Fun'(\DQ(\XX)^{\times n},\DQ(S))\to \Fun''(\DP(\XX)^{\times n},\DQ(S))
\]
is a categorical equivalence.

Let $\Fun^{\diamond}(\DQ(\XX)^{\times n},\DQ(S))$
be the full subcategory of $\Fun'(\DQ(\XX)^{\times n},\DQ(S))$,
spanned by functors which are compatible with
full subcategories $\DP(\XX)^{\times n}$ and $\DP(S)$.
Then the restriction functor
\[
\Fun^{\diamond}(\DQ(\XX)^{\times n},\DQ(S))\to \Fun''(\DP(\XX)^{\times n},\DP(S))
\]
is a categorical equivalences.
\end{Lemma}

\Proof
We first consider the case of $n=1$.
By \cite[5.3.5.10]{HTT} we have 
$\Fun_{\textup{cont}}(\CCC,\DDD)\simeq \Fun(\CCC_{\circ},\DDD)$,
where $\Fun_{\textup{cont}}(\CCC,\DDD)$ is the full subcategory of
$\Fun(\CCC,\DDD)$ spanned by functors that preserve filtered colimits.
Using \cite[5.3.5.15]{HTT} we see that if
$\Phi\in \Fun_{\textup{cont}}(\CCC,\DDD)$ is a left Kan extension
of $\phi\in\Fun''(\CCC_{\circ},\DDD)$, then
$\Phi$ preserves (co)kernels, that is, $\Phi$ is colimit-preserving.
Since the inclusion
$\mathcal{C}_{\circ}\to\CCC$ is exact, thus
we have an equivalence
$\Fun'(\CCC,\DDD)\simeq \Fun''(\CCC_{\circ},\DDD)$.

Now suppose that our claim holds
in the case $n=l$.
We will show that the map $\Fun'(\CCC^{\times (l+1)},\DDD)\to \Fun''(\CCC_{\circ}^{\times (l+1)},\DDD)$ is fully faithful.
Let $\mathcal{P}r^L$ be the $\infty$-category of presentable categories
in which functors are left adjoints.
Then by \cite[6.3.1.14]{HA}, $\mathcal{P}r^L$
has a symmetric monoidal structure ($\otimes$ denotes the tensor operation).
Then we have equivalences
\[
\Fun'(\CCC^{\times (l+1)},\DDD)\simeq \Fun^L(\CCC^{\otimes (l+1)},\DDD)\simeq \Fun^L(\CCC,\Fun^L(\CCC^{\otimes l},\DDD))\simeq\Fun^L(\CCC,\Fun''(\CCC_{\circ}^{\times l},\DDD)).
\]
The notation $\Fun^L(\ ,\ )$ indicates the full subcategory spanned by
left adjoints and $\CCC^{\otimes (l+1)}$ is the $(l+1)$-times product $\CCC\otimes \cdots \otimes \CCC$.
The above first equivalence follows from
the definition of $\CCC\otimes\CCC$ (see \cite[6.3.1]{HA}).
The second equivalence follows from the closed monoidal structure
of $\mathcal{P}r^L$ (cf. \cite[5.5.3.9]{HTT}).
The third one follows from the case of $n=l$.
Then by \cite[5.5.3.10]{HTT} again, we have
$\Fun_{\textup{cont}}(\CCC,\Fun(\CCC_{\circ}^{\times l},\DDD))\simeq\Fun(\CCC_{\circ},\Fun(\CCC_{\circ}^{\times l},\DDD))$ and they contain
$\Fun^L(\CCC,\Fun''(\CCC_{\circ}^{\times l},\DDD))$ as full subcategories.
Thus we have a functor
$\Fun'(\CCC^{\times (l+1)},\DDD)\to \Fun(\CCC_{\circ}^{\times (l+1)},\DDD)$
which is fully faithful.
Next we show that 
$\Fun'(\CCC^{\times (l+1)},\DDD)\to \Fun''(\CCC_{\circ}^{\times (l+1)},\DDD)$ is essentially surjective.
Since $\CCC_{\circ}\to \CCC$ is exact, the essential image of
$\Fun'(\CCC^{\times (l+1)},\DDD)\to \Fun(\CCC_{\circ}^{\times (l+1)},\DDD)$
lies in $\Fun''(\CCC_{\circ}^{\times (l+1)},\DDD)$.
Let $\CCC^{\times (l+1)}\to \DDD$
be a left Kan extension of $\CCC_{\circ}^{\times (l+1)}\to \DDD$,
which preserves filtered colimits separately in
each variable of $\CCC\times\cdots\times\CCC$.
To prove that
$\Fun'(\CCC^{\times (l+1)},\DDD)\to \Fun''(\CCC_{\circ}^{\times (l+1)},\DDD)$
is essentially surjective, it is enough to observe
that if $\CCC_{\circ}^{\times (l+1)}\to \DDD$ preserves finite colimits
separately in each variable, then the Kan extension
$\CCC^{\times (l+1)}\to \DDD$ preserves colimits separately in each variable.
It suffices to check it in each variable separately. Thus it follows from
the case of $n=1$.

Now the latter assertion is clear because $\CCC_{\circ}\to \CCC$ and
$\DDD_{\circ}\to \DDD$ are exact.
\QED

\begin{Lemma}
\label{looping}
Let $\mathcal{E}$ be a stable $\infty$-category.
The inclusion $\DP^{\le 0}(\XX)\to \DP(\XX)$ induces a fully faithful functor
\[
\Fun''(\DP(\XX),\mathcal{E})\to \Fun (\DP^{\le0}(\XX),\mathcal{E}).
\]
\end{Lemma}

\Proof
It follows along the same lines as the proof of \cite[1.3.3.11]{HA}.
\QED

Let $\Map_{\dagger}(\DP(\XX),\DP(S))$
be the full subcategory spanned by exact functors $\Phi: \DP(\XX) \rightarrow \DP(S)$
such that $\Phi(\DV(\XX))$ lies in $\DV(S)$.
By Lemma~\ref{leftkan} and~\ref{looping} together with
Proposition~\ref{cofinality}, we obtain:

\begin{Corollary}
\label{extension}
Suppose that $\XX$ has cohomological dimension zero.
Then the functor
\[
\Map_{\dagger}(\DP(\XX),\DP(S)) \to \Map(\DV(\XX),\DV(S))
\]
is a fully faithful functor.
\end{Corollary}

\begin{Lemma}
\label{extensionstep3}
Let $\mathcal{E}$ be a stable $\infty$-category.
Then the natural functor
\[
\Fun''(\DP(\XX)^{\times n},\mathcal{E})\to \Fun((\DP(\XX)^{\le0})^{\times n},\mathcal{E})
\]
is fully faithful.
\end{Lemma}

\Proof
Let $\CCC^{\le0}_{\circ}=\DP^{\le0}(\XX)$.
The case of $n=1$ follows from Lemma~\ref{looping}.
Now suppose that the case of $n=l$ holds.
There are natural fully faithful functors
\[
\Fun''(\CCC_{\circ}^{\times (l+1)},\mathcal{E})\to \Fun''(\CCC_{\circ},\Fun(\CCC_{\circ}^{\times l},\mathcal{E}))\to \Fun(\CCC^{\le0}_{\circ},\Fun(\CCC_{\circ}^{\times l},\mathcal{E})).
\]
The second functor is fully faithful by Lemma~\ref{looping} and the fact that $\Fun(\CCC_{\circ}^{\times l},\mathcal{E})$
is stable.
The essential image of $\Fun''(\CCC_{\circ}^{\times (l+1)},\mathcal{E})$
lies in $\Fun(\CCC^{\le0}_{\circ},\Fun''(\CCC_{\circ}^{\times l},\mathcal{E}
)).$
In addition by the case of $n=l$
we have a fully faithful functor
\[
\Fun(\CCC^{\le0}_{\circ},\Fun''(\CCC_{\circ}^{\times l},\mathcal{E}))\to \Fun(\CCC_{\circ}^{\le0},\Fun((\CCC_{\circ}^{\le0})^{\times l},\mathcal{E})).
\]
Since $\Fun(\CCC_{\circ}^{\le0},\Fun((\CCC_{\circ}^{\le0})^{\times l},\DDD_{\circ}))\simeq\Fun((\CCC_{\circ}^{\le0})^{\times (l+1)},\DDD_{\circ})$, thus the case of $n=l+1$ follows.
\QED

\begin{Lemma}
\label{extensionstep5}
Suppose that $\XX$ has cohomological dimension zero.
Then the restriction
\[
\Fun''(\DP(\XX)^{\times n},\DP(S))\to \Fun(\DV(\XX)^{\times n},\DP(S))
\]
is a fully faithful functor.
\end{Lemma}

\Proof
Let $\mathcal{E}$ be a stable presentable
$\infty$-category.
We may replace $\DP(S)$ by $\mathcal{E}$ (consider $\textup{Ind}\DP(S)$).
We first consider the case of $n=1$.
By Lemma~\ref{cofinality} (and its proof),
for any $P\in \CCC_{\circ}^{\le0}$, $P$ is a colimit of the natural diagram $(\CCC_{v})_{/P}\to \CCC$.
Since any object $F$ in the full subcategory $\Fun''(\CCC_{\circ},\mathcal{E})\subset\Fun (\CCC_{\circ}^{\le0},\mathcal{E})$ (cf. Lemma~\ref{extensionstep3})
extends to a colimit-preserving functor $\CCC \to \mathcal{E}$ by Lemma~\ref{extensionstep1},
$F$ is a left Kan extension of $F|_{\CCC_v}$.
Thus we have a fully faithful embedding $\Fun''(\CCC_{\circ},\mathcal{E})\subset\Fun(\CCC_v,\mathcal{E})$ induced by the inclusion
$\CCC_v\to \CCC_{\circ}$.

Next suppose that the case of $n=l$ holds.
We have fully faithful functors
\[
\Fun''(\CCC_{\circ}^{\times (l+1)},\mathcal{E})\to \Fun''(\CCC_{\circ},\Fun(\CCC_{\circ}^{\times l},\mathcal{E}))\to \Fun(\CCC_{\circ}^{\le0},\Fun(\CCC_{\circ}^{\times l},\mathcal{E})).
\]
By the observation in the case of $n=1$
(note that $\Fun(\CCC_{\circ}^{\times l},\mathcal{E})$ is stable
and presentable), we have a fully faithful
embedding
\[
\Fun''(\CCC_{\circ},\Fun(\CCC_{\circ}^{\times l},\mathcal{E}))\subset 
\Fun(\CCC_{v},\Fun(\CCC_{\circ}^{\times l},\mathcal{E})).
\]
Note that if a functor $F:\mathcal{C}_v\to \Fun(\CCC_{\circ}^{\times l},\mathcal{E})$
is in
the essential image of $\Fun''(\CCC_{\circ}^{\times (l+1)},\mathcal{E})$
then $F(\mathcal{C}_v)$ maps to $\Fun''(\CCC_{\circ}^{\times l},\mathcal{E})$.
Using the case of $n=l$ we also have a fully faithful embedding
$\Fun(\CCC_v,\Fun''(\CCC_{\circ}^{\times l},\mathcal{E}))\subset \Fun(\CCC_v^{\times (l+1)},\mathcal{E})$. This completes the proof.
\QED

\begin{Lemma}
\label{openreduction}
Let $j:\UU\to \XX$ be a quasi-compact open immersion.
Then the restriction functor $j^*:\DQ(\XX)\to \DQ(\UU)$
and $j_*:\DQ(\UU)\to \DQ(\XX)$
induce a localization (cf. \cite[5.2.7.2]{HTT})
\[
j^*:\DQ(\XX)\rightleftarrows \DQ(\UU):j_*.
\]
\end{Lemma}

\Proof
Let $T$ be a collection of morphisms $F\to F'\in \Fun (\Delta^1,\CCC)$
which are quasi-isomorphisms on $\UU$.
We let $T^{-1}\CCC$ be the full subcategory of $\CCC$
spanned by $T$-local objects, that is, objects $F\in \CCC$ such that
$\Map_{\CCC}(E',F)\to \Map_{\CCC}(E,F)$ is a homotopy equivalence
for any $E\to E'\in T$.
We claim  that $\CCC'$ is equivalent to $T^{-1}\CCC$.
More precisely, we will observe that
$j_*:\CCC'\to \CCC$ factors through $T^{-1}\CCC$
and it is a homotopy inverse of $j^*:T^{-1}\CCC\to \CCC'$.
To see that $j_*E$ is a $T$-local object
for any $E\in \CCC'$,
it suffices to show that, if $F'\to F$ in $\CCC$
is a quasi-isomorphism
on $\UU$, then the
induced functor $\Map_{\CCC}(F,j_*E)\to \Map_{\CCC}(F',j_*E)$
is a weak homotopy equivalence.
This equivalence follows from weak homotopy equivalences
\[
\Map_{\CCC}(F,j_*E)\simeq \Map_{\CCC'}(j^*F,E),\ 
\Map_{\CCC}(F',j_*E)\simeq \Map_{\CCC'}(j^*F',E)\] induced by the adjunction,
and $\Map_{\CCC'}(j^*F,E)\simeq
\Map_{\CCC'}(j^*F',E)$.
Hence $j_*E$ is a $T$-local object.
It remains to show that $j_*:\CCC'\to T^{-1}\CCC$
is a homotopy inverse of $j^*:T^{-1}\CCC\to \CCC'$.
For any $E\in \CCC'$ the adjoint map $j^*j_*E\to E$ is a quasi-isomorphism
($j$ is an open immersion).
For any $F\in T^{-1}\CCC$,
the adjoint map $F\to j_*j^*F$ is a quasi-isomorphism on $\UU$ (this means
that $F\to j_*j^*F$ is an equivalence in $T^{-1}\CCC$).
Thus $T^{-1}\CCC\simeq\CCC'$.

\QED


\section{Symmetric monoidal functors and Derived Morita theory}\label{sec:4}


Let $\stk X$ be an algebraic stack over a field $k$
and let $S$ be a scheme over $k$.
Let $\mfn: \dqc(\stk X) \rightarrow \dqc(S)$ be a $k$-linear symmetric monoidal functor which preserves small colimits.
We first give a condition under which $\mfn$ preserves vector bundles,
i.e., $\mfn(\dvb(\stk X))$ lies in $\dvb(S)$.

Let us recall the notions of \textit{integral functors} and their \textit{integral kernels}.
An object $P\in \dqc(\stk X \times S)$ gives rise to an exact functor $\mfn_{P}:=\ps2(\pl1(-)\otimes P)$
where $\mathrm{pr_1}$ and $\mathrm{pr_2}$ denote the natural projections from $\stk X \times S$ to $\stk X$ and $S$ respectively.
The functor $\mfn_{\cpx P}$ is  called \textit{integral functor}
and $P$ is called an \textit{integral kernel}, or simply \textit{kernel} of $\mfn_{P}$.
To avoid unnecessary confusion
we often denote by $\lotimes$ the derived tensor operation
 and denote by $\otimes$ the ordinary tensor operation.
Similarly, $\mathbb{R}(\bullet)_*$ means the derived pushforward,
whereas $(\bullet)_*$ indicates the ordinary pushforward.
Moreover, to emphasize that an object is a cochain complex
we often write $P^{\bullet}, Q^{\bullet},\ldots$ for cochain complexes.
We write $\textup{D}_{\textup{qcoh}}(\bullet)$ for the homotopy category
$\textup{h}\DQ(\bullet)$.
When we emphasize that
$\DQ(\bullet)$
(resp. $\textup{D}_{\textup{qcoh}}(\bullet)$)
equips with
the natural symmetric monoidal
structure, we then denote by
$\DQ^\otimes(\bullet)$
(resp. $\textup{D}_{\textup{qcoh}}^\otimes(\bullet)$).
If there is no danger of confusion,
we sometimes omit the subscript $\otimes$.

\begin{Proposition}\label{lbvb}
Let $\stk X$ be an algebraic stack over $k$.
Suppose that the cohomological dimension of $\stk X$ is finite,
i.e., there exists an integer $d$ such that for any quasi-coherent $\OO_{\stk X}$-module $F$ and $q>d$, we have $\tH^q(\stk X, F)=0$.
Let $S$ be a scheme over $k$.
Let $\mfn : \tdqc^\otimes(\stk X) \rightarrow \tdqc^\otimes(S)$ be a symmetric monoidal functor
whose underlying functor is an integral functor
induced by a bounded kernel $\cpx P \in \tdqc^\mathrm{b}(\stk X \times_k S)$.
Then $\mfn$ preserves vector bundles.
\end{Proposition}

Before the proof of this proposition, let us recall the notion of \textit{dualizable} objects in a symmetric monoidal category.
Let $(\mathcal C, \otimes, \bu)$ be a (ordinary) symmetric monoidal category.
An object $M$ in $\mathcal C$ is called \textit{dualizable}
if there exist an object $M^* \in \mathcal C$ and morphisms $\eta : \bu \rightarrow M \otimes M^*$
and $\epsilon : M^* \otimes M \rightarrow \bu$ satisfying the following conditions:
\begin{itemize}
\item The composite
$M \xrightarrow{\eta \otimes \id_M} M \otimes M^* \otimes M \xrightarrow{\id_M \otimes \epsilon} M$ is the identity map,
\item The composite
$M^* \xrightarrow{\id_{M^*} \otimes \eta} M^* \otimes M \otimes M^* \xrightarrow{\epsilon \otimes \id_{M^*}} M^*$ is the identity map.
\end{itemize}
The object $M^*$ is called a \textit{dual} of $M$.
If $M^*$ exists, it is unique up to isomorphism.
If $\Psi: \mathcal C \rightarrow \mathcal C'$ is a symmetric monoidal functor and $M$ is a dualizable object of $\mathcal C$,
$\Psi(M)$ is also dualizable and $\Psi(M)^* \simeq \Psi(M^*)$.
In the case that $\mathcal C$ is a category of quasi-coherent complexes,
any perfect complex $\cpx E$ is dualizable and its dual is isomorphic to the (usual) derived dual $\RR\mathcal Hom(\cpx E, \OO)$.
Therefore, for any symmetric monoidal functor $\mfn: \tdqc(\stk X) \rightarrow \tdqc(S)$ and any perfect complex $\cpx E \in \tdqc(\stk X)$,
$\mfn(\RR\mathcal Hom(\cpx E, \OO_{\stk X}))$ is isomorphic to $\RR\mathcal Hom(\mfn(\cpx E), \OO_S)$.

\begin{Remark}
Let $\XX$ and $S$ be algebraic stacks.
Then any functor $\Phi:\DQ^\otimes(\XX)\to \DQ^\otimes(S)$ which is symmetric monoidal,
preserves dualizable objects.
According to \cite[Proposition 3.6]{BFN}
dualizable objects and perfect complexes coincide in $\DQ^\otimes(\XX)$.
Also, dualizable objects and perfect complexes
coincide in $\DQ^\otimes(S)$.
Consequently, any symmetric monoidal functor $\Phi$ preserves
perfect complexes.
\end{Remark}

\renewcommand{\proofname}{Proof of Proposition \ref{lbvb}.}

\begin{proof}

We may and will assume that $S$ is affine.
Let $d$ be the cohomological dimension of $\stk X$ and $m := \max \{\, p \mid \tH^p(\cpx P) \neq 0 \,\}$.
To prove this proposition, we first claim that $\tH^q(\mfn(E)) = 0$
for any vector bundles $E$ on $\stk X$ and $q > m+d$.
The category of quasi-coherent $\OO_{\stk X\times S}$-modules has enough injective objects.
(For the readers' convenience, we give an outline of the proof here.
Let $F$ be a quasi-coherent $\OO_{\stk X\times S}$-module
and let $p:U\rightarrow \stk X\times S$ be a smooth surjective map where $U$ is an affine scheme.
Take an injective quasi-coherent $\OO_U$-module $I$ which contains $p^*F$.
Since $p_*I$ is an injective $\OO_{\stk X \times S}$-module,
it is sufficient to check that the natural maps $F\rightarrow p_*p^*F$ and $p_*p^*F\rightarrow p_*I$ are injective.
The first follows from the fact that $p$ is faithfully flat and affine.
The second is clear.)
Hence there exists a bounded below complex of injective quasi-coherent $\OO_{\stk X\times S}$-modules $\cpx I^{\bullet}$
which is quasi-isomorphic to $\cpx P$.
Since $\pl1E$ is a vector bundle,
$\pl1E \otimes \cpx I^{\bullet}$ is quasi-isomorphic to $\pl1E \lotimes \cpx P$
and $\pl1E \otimes I^l$ is an injective quasi-coherent $\OO_{\stk X \times S}$-module for any $l\in\ZZ$.
Thus, we have
\begin{equation}\label{phiE}
\tH^q(\mfn(E))
=\tH^q(\RR\ps2(\pl1E\lotimes\cpx P))
\simeq \tH^q(\ps2(\pl1E\otimes\cpx I^{\bullet})).
\end{equation}
On the other hand, since $\tH^l(\cpx I^{\bullet})\simeq\tH^l(\cpx P) = 0$ for any $l>m$ and $\pl1E$ is flat,
we have $\tH^l(\pl1E \otimes \cpx I^{\bullet}) \simeq \pl1E \otimes \tH^l(\cpx I^{\bullet}) = 0$ for any $l>m$.
Hence the complex
\begin{equation}\label{ZEres}
0 \rightarrow \pl1E \otimes Z^m \rightarrow \pl1E \otimes I^m
\rightarrow \pl1E \otimes I^{m+1} \rightarrow \pl1E \otimes I^{m+2} \rightarrow \cdots
\end{equation}
is exact, where $Z^m$ is $\ker(I^m \rightarrow I^{m+1})$.
Moreover, since $\pl1E \otimes I^l $ is injective for any $l\in\ZZ$,
(\ref{ZEres}) gives an injective resolution of $\pl1E \otimes Z^m$.
Thus we have
\begin{equation}\label{ps2ZE}
\tH^q(\ps2(\pl1E \otimes \cpx I^{\bullet}))
\simeq \tH^{q-m}(\RR\ps2(\pl1E \otimes Z^m)).
\end{equation}
Since $q-m>d$, we have $\tH^q(\ps2(\pl1E \otimes \cpx I^{\bullet})) \simeq \tH^{q-m}(\RR\ps2(\pl1E \otimes Z^m)) = 0$.
Therefore we obtain $\tH^q(\mfn(E)) = 0$ by (\ref{phiE}) and (\ref{ps2ZE}).

We then show that $\tH^q(\mfn(E)) = 0$ for any vector bundle $E$ on $\stk X$ and $q > 0$.
If $\mfn(E) = 0$, we have nothing to prove, so we assume that $\mfn(E)\neq 0$. 
Let $l$ be the integer $\max\{\, q \mid \tH^q(\mfn(E)) \neq 0\,\}$.
In general, if $\cpx F$ and $\cpx G$ are objects in $\tdpf(S)$
such that $\tH^i(\cpx F) \simeq \tH^j(\cpx G) = 0$ for any $i>s$ and $j>t$,
then we have
$\tH^k(\cpx F \lotimes \cpx G) = 0$
for any $k>s+t$, and
\begin{equation}\label{abc}
\tH^{s+t}(\cpx F\lotimes\cpx G) \simeq \tH^s(\cpx F)\otimes\tH^t(\cpx G).
\end{equation}
(To see this, take a complex $\cpx A$ (resp. $\cpx B$) which is quasi-isomorphic to $\cpx F$ (resp. $\cpx G$)
such that $A^i$ is a \textit{flat} $\OO_S$-module for any $i\in\ZZ$ and $A^i=0$ for any $i>s$ (resp. $B^j=0$ for any $j>t$)
and compute the cohomologies of the total complex of the double complex $\cpx A \otimes \cpx B$,
which is quasi-isomorphic to $\cpx F \lotimes \cpx G$.)
Hence for any positive integer $n$, we have $\tH^{nl}(\mfn(E)^{\lotimes n}) \simeq \tH^l(\mfn(E))^{\otimes n}$
($\otimes n$ represents the $n$-times product). 
On the other hand, since $\mfn$ is symmetric monoidal, we have $\mfn(E)^{\lotimes n} \simeq \mfn(E^{\lotimes n})$.
Hence we have $\tH^{nl}(\mfn(E^{\lotimes n})) \simeq \tH^l(\mfn(E))^{\otimes n}$.
Since $\tH^l(\mfn(E))$ is a non-zero quasi-coherent sheaf of finite type by \cite[2.2.3]{TT},
it follows that $\tH^l(\mfn(E))^{\otimes n} \neq 0$.
Indeed, if $\tH^l(\mfn(E))^{\otimes n} = 0$,
then $(\tH^l(\mfn(E)) \otimes k(s))^{\otimes n} \simeq \tH^l(\mfn(E))^{\otimes n} \otimes k(s) = 0$ for any point $s \in S$
where $k(s)$ denotes the residue field of $s$.
This implies that $\tH^l(\mfn(E)) \otimes k(s) = 0$ since $\tH^l(\mfn(E)) \otimes k(s)$ is a $k(s)$-vector space.
Hence the stalk $\tH^l(\mfn(E))_s$ is zero by Nakayama's lemma and so $\tH^l(\mfn(E)) = 0$.
Therefore $\tH^{nl}(\mfn(E^{\lotimes n})) \simeq \tH^l(\mfn(E))^{\otimes n} \neq 0$.
We have to show that $l$ is not positive.
If $l$ is positive, there exists a positive integer $n$ such that $nl > m + d$.
In addition, since $E^{\lotimes n}$ is a locally free sheaf, $\tH^q(\mfn(E^{\lotimes n})) = 0$ for any $q > m+d$.
It gives rise to a contradiction.

Next, we show that $\tH^{-q}(\mfn(E)) = 0$ for any $q > 0$.
We have equivalences $\mfn(E) \simeq \mfn(E^{**}) \simeq \mfn(E^*)^* \simeq \RR\mathcal Hom(\mfn(E^*), \OO_S)$.
The second equivalence follows from the fact that $\mfn$ is symmetric monoidal.
On the other hand, since $E^*$ is a locally free sheaf, we have $\tH^q(\mfn(E^*)) = 0$ for any $q > 0$.
Hence $\mfn(E^*)$ is quasi-isomorphic to a complex $\cpx M$ such that $M^q=0$ for any $q>0$ and $M^q$ is a free module for any $q$
since $S$ is affine.
Thus we have
\begin{equation}\label{dualrep}
\tH^{-q}(\RR\mathcal Hom(\mfn(E^*), \OO_S)) \simeq \tH^{-q}(\mathcal Hom(\cpx M, \OO_S)) = 0.
\end{equation}
Therefore we have $\tH^{-q}(\mfn(E)) = 0$ for any $q > 0$ by (\ref{dualrep}).

It remains to prove that $\mfn(E) \simeq \tH^0(\mfn(E))$ is a vector bundle.
Since $\tH^0(\mfn(E))$ is quasi-coherent of finite type, it is enough to show that $\tH^0(\mfn(E))$ is flat.
To see this, it is enough to show that $\mathcal Tor_1^{\OO_S}(\tH^0(\mfn(E)), N)=0$
for any quasi-coherent $\OO_S$-module $N$. 
We have
\begin{eqnarray*}
\mathcal Tor_1^{\OO_S}(\tH^0(\mfn(E)), N)
&\simeq& \tH^{-1}(\tH^0(\mfn(E)) \lotimes N)\\
&\simeq& \tH^{-1}(\mfn(E) \lotimes N)\\
&\simeq& \tH^{-1}(\mathbb R\mathcal Hom (\mfn(E^*), N))\\
&\simeq& \tH^{-1}(\mathbb R\mathcal Hom (\tH^0(\mfn(E^*)), N)) = 0. 
\end{eqnarray*}
Therefore $\tH^0(\mfn(E))$ is flat and it is a locally free sheaf.
\end{proof}

\begin{Remark}
We will apply Proposition~\ref{lbvb} only to schemes $\XX$
in this paper.
\end{Remark}

\begin{Remark}\label{uba}
By the argument in the proof of Proposition \ref{lbvb},
we see the following: if $\mfn(\dvb(\stk X))$ is uniformly bounded above
(i.e., there exists an integer $a$ such that for any vector bundle $E$,
$\tH^l(\mfn(E))$ is zero for any $l>a$), then
$\mfn(E)$ is a vector bundle.
\end{Remark}

\renewcommand{\proofname}{Proof.}

\renewcommand{\proofname}{Proof.}

Let us recall one of key ingredients: {\it derived
Morita theory}
due to To\"en, which was further generalized by Ben-Zvi, Francis and Nadlar
(see \cite[Theorem 8.9]{To}, \cite[Corollary 4.10]{BFN}).
We here recall the form
which we can apply to our situation. 
Let $\stk X$ be a perfect stack over $k$.
Then there is a natural functor
\begin{equation}\label{lbFM}
\dqc(\stk X\times_k S) \rightarrow \Fun_{\mathcal{M}od_{\DQ(k)}(\mathcal{P}r^L)}(\dqc(\stk X), \dqc(S));\ \ P \mapsto \mfn_{P},
\end{equation}
where $\mathcal{M}od_{\DQ(k)}(\mathcal{P}r^L)$ is the $\infty$-category
of left $\DQ(k)$-modules in $\mathcal{P}r^L$
(see Section \ref{subsec:2-2}). Here $\stk X\times_kS$ is the fiber product in the category
of ordinary stacks, but it coincides with the fiber product of derived stacks
since $k$ is a field.

\begin{Theorem}[\cite{To}, \cite{BFN}]\label{lbToen}
Let $\stk X$ be a perfect algebraic stack over $k$.
Then $($\ref{lbFM}$)$ gives a categorical equivalence.
\end{Theorem}

\begin{Proposition}\label{lbker}
Let $X$ be a noetherian scheme endowed with a very ample invertible sheaf over $k$
and let $S$ be a scheme over $k$.
Let $\Phi:\textup{D}^\otimes_{\textup{qcoh}}(X)\to \textup{D}^\otimes_{\textup{qcoh}}(S)$ be a symmetric monoidal functor
whose underlying functor is an integral functor induced by
an integral kernel $\cpx P\in \tdqc(X\times_kS)$.
Then $\cpx P$ is a sheaf, that is, $\tH^l(\cpx P)=0$ for any $l\neq 0$.
\end{Proposition}

In the proof of this proposition,
we consider derived pushforwards of \textit{unbounded} complexes,
so let us recall the notion of \textit{K-injective} complexes (cf. \cite{Spal}).
A (unbounded) complex $\cpx A$ in an abelian category $\mathcal A$ is called \textit{K-injective}
if, for any acyclic complex $\cpx B$ in $\mathcal A$,
the complex $\hhom_\mathcal A^\bullet(\cpx B, \cpx A)$ is acyclic.
If $\mathcal A$ is the category of quasi-coherent sheaves on a scheme,
any complex in $\mathcal A$ is quasi-isomorphic to a K-injective complex.
For any morphism $f$ of schemes, the derived pushforward $\RR f_*\cpx E$ of a complex $\cpx E$ of quasi-coherent modules is quasi-isomorphic
to the non-derived pushforward $f_*\cpx I$ of a K-injective complex $\cpx I$ which is quasi-isomorphic to $\cpx E$. 

\renewcommand{\proofname}{Proof of Proposition \ref{lbker}.}

\begin{proof}
This problem is local on $S$, we may assume that $S$ is a connected affine
scheme.
Taking a K-injective resolution, we may assume that $\cpx P$ is K-injective.
For any $l \in \ZZ$, let $d^l$ be the differential map $P^l \rightarrow P^{l+1}$
and $\alpha^l: \ker d^l \rightarrow \tH^l(\cpx P)$ be the natural surjection.
To prove this proposition, it is enough to show that $\alpha^l = 0$ for any integer $l \neq 0$ (since $\alpha^l$ is surjective).

Let $\OO_X(1)$ be a very ample invertible sheaf on $X$
and let $\cpx Q(m)$ denote $\cpx Q \otimes \pl1 \OO_X(m)$ for any (unbounded) complex $\cpx Q$ of quasi-coherent $\OO_{X\times S}$-modules
on $X \times S$. Fix $l\ne 0$.
Now suppose that $\alpha^l \neq 0$.
Then there exist $f \in \Gamma(X \times S, \OO_{X \times S}(1))$ and $\phi \in \Gamma((X \times S)_f, \ker d^l)$
such that $\alpha^l|_{(X \times S)_f}(\phi) \neq 0$, where $(X \times S)_f$ denotes the affine open subscheme of $X \times S$ where $f$ does not vanish.
For any sufficiently large $n \in \ZZ$, $f^n\phi$ lies in $\Gamma(X \times S, \ker d^l(n))$
and thus it follows that $\Gamma(X \times S, \alpha^l(n)) \neq 0$,
where $\alpha^l(n)$ denotes $\alpha^l \otimes \id_{\OO_{X \times S}(n)} : \ker d^l(n) \rightarrow \tH^l(\cpx P)(n)$.
Hence, to see that $\alpha^l = 0$,
it is enough to show that the induced morphism
$\Gamma(X \times S, \alpha^l(N)): \Gamma(X \times S, \ker d^l(N)) \rightarrow \Gamma(X \times S, \tH^l(\cpx P(N))) $ is zero
for any sufficiently large $N \in \ZZ$.
Since $S$ is affine, this is equivalent to showing that the induced morphism
$\ps2(\alpha^l(N)): \ps2(\ker d^l(N)) \rightarrow \ps2\tH^l(\cpx P(N))$ is zero
where $\ps2$ denotes the \textit{non}-derived pushforward.
Applying $\ps2$ to the complex
\[
\cpx P(N) : \cdots \rightarrow P^{l-1}(N) \xrightarrow{d^{l-1}(N)} P^l(N) \xrightarrow{d^l(N)} P^{l+1}(N) \rightarrow \cdots
\]
we obtain a complex
\[
\ps2(\cpx P(N)):
\cdots \rightarrow \ps2P^{l-1}(N) \xrightarrow{\ps2(d^{l-1}(N))} \ps2P^l(N) \xrightarrow{\ps2(d^l(N))} \ps2P^{l+1}(N) \rightarrow \cdots.
\]
From these complexes we have the following commutative diagram:
\[
\xymatrix{
\ker (\ps2 d^l(N)) \ar[rr] && \tH^l(\ps2(\cpx P(N)))\ar[d] \\
\ps2 (\ker d^l(N)) \ar[rr]^{\ps2(\alpha^l(N))} \ar[u]_{\simeq} && \ps2 \tH^l(\cpx P(N)). \\
}
\]
Hence, to show that $\ps2(\alpha^l(N)) = 0$,
it is enough to show that $\tH^l(\ps2(\cpx P(N))) = 0$ for any sufficiently large integer $N$.
Since $P$ is a K-injective complex and $\pl1\OO_X(N)$ is invertible,
$\cpx P(N)$ is also a K-injective complex
and hence
$\ps2(\cpx P(N))$
is quasi-isomorphic to $\RR\ps2(\pl1\OO_X(N) \lotimes \cpx P)$.
Moreover, since $\cpx P$ is an integral kernel of $\mfn$,
we have $\mfn(\OO_X(N)) \simeq \RR\ps2(\pl1\OO_X(N) \lotimes \cpx P)$.
Hence $\mfn(\OO_X(N))$ is quasi-isomorphic to the complex $\ps2(\cpx P(N))$.
Thus, to show that $\tH^l(\ps2(\cpx P(N))) = 0$, it will suffice to show that $\tH^l(\mfn(\OO_X(N))) = 0$.
By \cite[Theorem 2.3]{Yek} and the connectedness of $S$, for any two objects $\cpx F_1$, $\cpx F_2 \in \tdpf(S)$
such that $\cpx F_1 \lotimes \cpx F_2 \simeq \OO_S$, there exist an invertible sheaf $L$ on $S$ and $m \in \ZZ$
such that $\cpx F_1 \simeq L[m]$.
Since $\mfn$ preserves $\lotimes$ and structure sheaves, 
there exist an invertible sheaf $L$ on $S$ and $m \in \ZZ$ such that $\mfn(\OO_X(1)) \simeq L[m]$
and we have $\mfn(\OO_X(N)) \simeq L^{\lotimes N}[Nm]$.
Since $\tH^l(L^{\lotimes N}[Nm]) = 0$ if $l \neq -Nm$,
thus $\tH^
l(\mfn(\OO_X(N))) = 0$ for a sufficiently large $N \in \ZZ$.
\end{proof}

\renewcommand{\proofname}{Proof.}

\begin{Corollary}
\label{schemepreservevect}
Let $\XX$ be a scheme that satisfies (i)
in Section \ref{subsec:2-3}. Let $S$ be a scheme over $k$.
Let $\Phi:\textup{D}_{\textup{qcoh}}^{\otimes}(\XX)\to \textup{D}_{\textup{qcoh}}^\otimes(S)$ be a symmetric
monoidal whose underlying functor is an integral functor
induced by an integral kernel in $\textup{D}_{\textup{qcoh}}(\XX\times_kS)$.
Then $\Phi$ preserves vector bundles.
\end{Corollary}

Next we consider the case (ii).


\begin{Proposition}\label{lbpvb}
Let $\stk X$ be an algebraic stack that satisfies $(ii)$ in Section \ref{subsec:2-3}.
Let $S$ be a scheme over $k$.
Let $\Phi:\textup{D}^\otimes_{\textup{qcoh}}(\XX)\to \textup{D}^\otimes_{\textup{qcoh}}(S)$ be a symmetric
monoidal functor.
Suppose that the underlying functor of $\Phi$ is an integral functor
induced by an integral kernel in $\textup{D}_{\textup{qcoh}}(\XX\times_kS)$.
We abusively denote the integral functor by the same symbol $\Phi:\dqc(\stk X)\to\dqc(S)$.
Then $\Phi$ preserves vector bundles.
\end{Proposition}

\begin{proof}
For simplicity of notation, in this proof we denote by $\otimes$ (resp. $f^*$)
the derived tensor operation
(resp. derived pullback functor).
We may suppose that $S$ is affine.

\textbf{Case 1.}
First we assume that $k$ is algebraically closed and $S$ is $\Spec k$.
We will show that there exists a closed point $\tilde{x}$ of $\stk X$ such that for any vector bundle $E$ on $\stk X$,
$\mfn(E)$ is determined by the restriction of $E$ to $\tilde{x}$.
Let $p: \stk X \rightarrow M$ denote the coarse moduli map.
Since $\mfn \circ p^*$ is the composite of an integral functor
and $p^*$, and $M$ satisfies the condition (i) in Section \ref{subsec:2-3}, thus by
Corollary~\ref{schemepreservevect}, Theorem~\ref{full} and Proposition~\ref{combining} (see Remark~\ref{inadequate})
there exists a morphism $x: S=\Spec k \rightarrow M$ such that $x^* \simeq \mfn \circ  p^*$. 
This morphism $x$ determines a closed point of $M$ which we denote by the same letter $x$.
Let $\co_{M,x}$ be the completion of the local ring $\OO_{M,x}$
and let $\co_{\stk X,p^{-1}(x)}$ be the completion of $\OO_{\stk X}$
with respect to the ideal $I$ of the closed substack $p^{-1}(x)$.
Since $\OO_{M,x}$ is noetherian, $\co_{M,x}$ is a flat $\OO_M$-module.
Thus we have $x^*\co_{M,x}\simeq k$ and $p^*\co_{M,x}\simeq\co_{\stk X,p^{-1}(x)}$.
Therefore we have
\[
\mfn(E) \simeq \mfn(E)\otimes x^*\co_{M,x}
\simeq \mfn(E)\otimes \mfn(p^*\co_{M,x})
\simeq \mfn(E\otimes p^*\co_{M,x})
\simeq \mfn(E\otimes \co_{\stk X,p^{-1}(x)}).
\]
This means that $\mfn(E)$ is determined by the pullback of $E$ to the stack $\stk X':=\stk X\times_M \Spec \co_{M,x}$.
Since $p$ is proper, by the Grothendieck's existence theorem for stacks \cite[Theorem 1.4]{Ols},
the category of coherent sheaves on $\stk X'$ is equivalent
to the category of compatible systems $\{(F'_n, \phi_n':F'_{n+1}/m^{n+1}F'_{n+1} \xrightarrow{\thicksim} F'_n)_{n\geq 0}\}$
of coherent sheaves on the reductions $\stk X'_n:=\stk X \times_{M} \Spec (\OO_{M,x}/m^{n+1})$
where $m$ is the maximal ideal of $\OO_{M,x}$
and $\phi'_n$ is an isomorphism of coherent sheaves.
Let $J \subset \OO_{\stk X}$ be the ideal of the closed substack $(\stk X'_0)_{\mathrm{red}}$
and $\stk X_n$ denote the closed substack defined by $J^{n+1}$.
Then the category of compatible systems of coherent sheaves on $\stk X'_n$ is equivalent to
the category of compatible systems of coherent sheaves on $\stk X_n$.
Therefore we can regard any vector bundle $E'$ on $\stk X'$
as a system $\{E_n\}_{n\geq 0}$ where $E_n$ is a vector bundle on $\stk X_n$
and $E_{n+1}$ is a flat deformation of $E_n$ to $\stk X_{n+1}$.
We will observe that this system $\{E_n\}_{n\geq 0}$ is determined by $E_0$.
According to the deformation theory of modules over a ringed topos
\cite[IV, Proposition 3.1.5]{Ill}, the set of isomorphism classes of flat deformations of $E_n$ to $\stk X_{n+1}$
is a torsor under $\mathrm{Ext}^1_{\OO_{\stk X_n}}(E_n, E_n\otimes_{\OO_{\stk X_n}} J^{n+1}/J^{n+2})$ (note that
Cartesian modules are stable under deformations) and we have
\begin{eqnarray*}
\mathrm{Ext}^1_{\OO_{\stk X_n}}(E_n, E_n\otimes_{\OO_{\stk X_n}} J^{n+1}/J^{n+2})
&\simeq& \tH^1(\stk X_n, \mathcal Hom_{\OO_{\stk X_n}}(E_n, E_n\otimes_{\OO_{\stk X_n}} J^{n+1}/J^{n+2}))\\
&\simeq& \tH^1(\stk X_0, \mathcal Hom_{\OO_{\stk X_0}}(E_0, E_0\otimes_{\OO_{\stk X_0}} J^{n+1}/J^{n+2})).
\end{eqnarray*}
Let $\tilde{x}: \Spec k \rightarrow \stk X$ be a point of $\stk X$ such that $p\circ \tilde{x}=x$.
Then $\stk X_0$ is isomorphic to the residual gerbe of $\tilde{x}$ over $k$
and this gerbe is isomorphic to the classifying stack $BG_{\tilde{x}}$,
where $G_{\tilde{x}}$ is the stabilizer group of $\tilde{x}$.
Since $\stk X$ is tame, $G_{\tilde{x}}$ is linearly reductive
and hence $\tH^1(\stk X_0, \mathcal Hom_{\OO_{\stk X_0}}(E_0, E_0\otimes_{\OO_{\stk X_0}} J^{n+1}/J^{n+2}))=0$.
Therefore the system $\{E_n\}$ is determined by $E_0$.
In addition, in our setting $G_{\tilde{x}}$ is finite.
Hence the number of finite dimensional irreducible representations of $G_{\tilde{x}}$ is finite
and any representation is completely reducible.
In other words, there exist vector bundles
$E_{01},\ldots, E_{0n}$ on $\stk X_0$
such that any vector bundle on $\stk X_0$ is isomorphic to a sheaf of the form $\bigoplus E_{0i}^{\oplus a_i}$.
By the deformation theory
and the Grothendieck's existence theorem,
for any $i$,
there exists an object $F_i$ in $\dqc(\stk X)$
which is a locally free $\co_{\stk X,p^{-1}(x)}$-module
and whose restriction to $\stk X_0$ is isomorphic to $E_{0i}$.
Thus for any vector bundle $E$ on $\stk X$,
$\mfn(E)$ is quasi-isomorphic to a complex of the form $\bigoplus\mfn(F_i)^{\oplus a_i}$.
If $a_i \neq 0$, the complex $\mfn(F_i)$ is bounded 
since $\mfn(E)$ is bounded for any vector bundle $E$.
Hence the family $\mfn(\dvb(\stk X))$ is uniformly bounded,
i.e., there exist integers $a \leq b$ such that for any vector bundle $E$,
$\tH^l(\mfn(E))=0$ if $l\not\in[a,b]$. 
Therefore, by Remark~\ref{uba},
it follows that $\mfn$ preserves vector bundles.

\textbf{Case 2.}
We then consider the case that $k$ is algebraically closed and $S$ is a general affine scheme over $k$.
We will prove that the family $\mfn(\dvb(\stk X))$ is uniformly bounded above (Remark \ref{uba}).
If this family is not uniformly bounded above, there exist a vector bundle $E$ on $\stk X$ such that
the integer $m = \mathrm{max}\{\, l \mid \tH^l(\mfn(E))\neq 0 \,\}$ is positive.
Since $\tH^m(\mfn(E))$ is finitely generated, by Nakayama's lemma,
there exist a field $K$ and a morphism $a: \Spec K \rightarrow S$ such that $\tH^m(a^*\mfn(E))$ is not zero.
Hence we may and will assume that $S=\Spec K$.
By Corollary~\ref{schemepreservevect}, Theorem~\ref{full} and Proposition~\ref{combining} there exists a morphism $f:S \rightarrow M$ such that $f^*\simeq \mfn\circ p^*$. 
Since $M$ is of finite type over $k$, there exist a $k$-subalgebra $R$ of $K$ of finite type and $g: T=\Spec R \rightarrow M$
such that $f=g\circ h$ where $h: S \rightarrow T$ is the morphism induced by the inclusion $R\subset K$.
We have obtained the following homotopy commutative diagram:
\begin{eqnarray*}
\begin{xymatrix}
{
\dqc(\stk X) \ar[r]^(0.4){\mfn} & \dqc(S=\Spec K) \\
\dqc(M) \ar[r]^(0.4){g^*} \ar[u]^{p^*} \ar[ru]^{f^*} & \dqc(T=\Spec R) \ar[u]^{h^*}. \\
}
\end{xymatrix}
\end{eqnarray*}
Let $\xi: \Spec k \rightarrow T$ be a closed point of $T$
and let $x$ denote the closed point $g\circ\xi:\Spec k \rightarrow M$.
By the similar argument as in Case 1,
there exist objects $F_1,\ldots,F_n$ in $\dqc(\stk X)$ such that $\mfn(F_i)$ is bounded above and
for any vector bundle $E$ on $\stk X$,
$E\otimes \co_{\stk X,p^{-1}(x)}$ is quasi-isomorphic to a complex of the form $\bigoplus_jF_j^{\oplus a_j}$. 
Hence $\mfn(E\otimes\co_{\stk X,p^{-1}(x)})$ is
quasi-isomorphic to a complex of the form $\bigoplus_j\mfn(F_j)^{\oplus a_j}$
and the family $\{\mfn(E\otimes\co_{\stk X,p^{-1}(x)})\}_{E\in\dvb(\stk X)}$ is uniformly bounded above.
On the other hand, we have
\[
\mfn(E\otimes\co_{\stk X,p^{-1}(x)})
\simeq \mfn(E)\otimes f^*\co_{M,x}
\]
Note that $f^*\co_{M,x}$ is not zero.
Hence the family $\mfn(\dvb(\stk X))$ is uniformly bounded above.

\textbf{Case 3.}
Here we consider the case of an arbitrary base field $k$.
Let $k\subset \ol{k}$ be an algebraic closure.
As in Case 2 we may and will assume that $S$ is $\Spec K$ where $K$ is a field.
For an algebraic stack $\mathcal{Y}$ over $k$
we will write $\ol{\mathcal{Y}}$ for $\mathcal{Y}\times_k\ol{k}$.
By \cite[Theorem 4.7.]{BFN}, $\dqc(\ol{\stk X})$
and $\dqc(\ol S)$ are naturally equivalent to
\begin{eqnarray*}
\dqc(\stk X) \otimes_{\dqc(\Spec k)} \dqc(\Spec \ol k)\ \ \ \textup{and}\ \ \dqc(S) \otimes_{\dqc(\Spec k)} \dqc(\Spec \ol k)
\end{eqnarray*}
respectively (see \cite{BFN} for the notation).
Thus we have $\ol{\Phi}=\Phi\otimes_{\DQ(k)}\DQ(\ol{k}):\DQ(\ol{\XX})\to \DQ(\ol{S})$.
Let $E$ be a vector bundle
on $\XX$.
It suffices to show that $\Phi(E)$ is (quasi-isomorphic to)
a locally free sheaf in
$\DQ(\ol{S})$, that is, $\ol{\Phi}(E)$ is a locally free sheaf.
To complete the proof, we will reduce this case to the Case 1 and 2.
Let $f:S\to M$ be a morphism such that $\Phi\circ p^*\simeq f^*$.
If $\ol{p}$ and $\ol{f}$ denote the base changes of the
coarse moduli map $p:\XX\to M$ and $f:S\to M$ respectively,
then $\ol{\Phi}\circ \ol{p}^*\simeq \ol{f}^*$
since external products of objects in $\DQ(M)$ and $\DQ(\ol{k})$
generate $\DQ(\ol{M})$ (see \cite[Section 4.2]{BFN}).
Let $F$ be a quasi-coherent sheaf on $\ol{\XX}$.
Since $F$ is an inductive colimit of
external products of objects in $\DQ(\XX)$ and $\DQ(\ol{k})$, thus it follows that
$\ol{\Phi}(E\otimes F)\simeq \ol{\Phi}(E)\otimes \ol{\Phi}(F)$.
Consequently, we can apply the arguments in Case 1 and 2 to $\ol{\Phi}$ and
complete the proof.
\end{proof}

\begin{Remark}
\label{inadequate}
The proof of Proposition~\ref{lbpvb} uses Proposition~\ref{combining}
of the case (i) of Section \ref{subsec:2-3}.
But the proof of Proposition~\ref{combining}
for the scheme case (i) does not need Proposition~\ref{lbpvb}.
\end{Remark}

As a consequence of this section, we here record the following corollary which
follows from Theorem~\ref{lbToen}, Corollary~\ref{schemepreservevect}, Proposition~\ref{lbpvb} and Theorem~\ref{full}.
See Section \ref{sec:5} for the notion of $\infty$-categorical
 symmetric monoidal functor, but readers who are not familiar with it might skip to the next section.

\begin{Corollary}
\label{summary4}
Let $\XX$ be an algebraic stack that satisfies either (i) or (ii)
in Section \ref{subsec:2-3}. Let $S$ be an affine scheme over $k$.
Let $\Phi:\DQ^{\otimes}(\XX)\to \DQ^\otimes(S)$ be a symmetric
monoidal colimit-preserving functor over $\DQ^{\otimes}(k)$.
Then $\Phi$ preserves vector bundles,
and there exist a $k$-morphism $f:S\to \XX$
and a monoidal natural transformation
$f^*|_{\DV^\otimes(\XX)}\simeq \Phi|_{\DV^\otimes(\XX)}$.
\end{Corollary}






\section{Derived Tannaka duality}\label{sec:5}
In this section using results of Section \ref{sec:3} and \ref{sec:4}
we prove main results of this paper Theorem~\ref{supermain},
Corollary~\ref{reconstruction} and Theorem~\ref{geomono}.

We here use the theory of symmetric monoidal $\infty$-categories
developed in \cite{HA}.
We refer to \cite{HA} for its generalities.
Let $\mathcal{F}\textup{in}_*$ be the category of marked finite sets
(our notation is slightly different from \cite{HA}).
Namely, objects are marked finite sets and
a morphism from $\langle n\rangle_*:=\{1<\cdots <n\}\sqcup\{*\}\to \langle m\rangle_*:=\{1<\cdots <m\}\sqcup\{*\}$
is a (not necessarily order-preserving) map of finite sets which preserves the distinguished points $*$.
Let $\alpha^{i,n}:\langle n\rangle_*\to \langle 1\rangle_*$ be
a map such that $\alpha^{i,n}(i)=1$ and $\alpha^{i,n}(j)=\ast$ if $i\neq j\in \langle n\rangle_*$.
Let $I:=\textup{N}(\mathcal{F}\textup{in}_*)$.
A symmetric monoidal category is a coCartesian fibration (cf. \cite[2.4]{HTT})
$p:\mathcal{M}^\otimes\to I$ such that
for any $n\ge0$, $\alpha^{1,n}\ldots \alpha^{n,n}$ induce
an equivalence $\mathcal{M}^\otimes_{n}\to (\mathcal{M}^\otimes_1)^{\times n}$
where $\mathcal{M}^\otimes_{n}$ and $\mathcal{M}^\otimes_{1}$
are fibers of $p$ over $\langle n\rangle_*$ and $\langle 1\rangle_*$
respectively.
A symmetric monoidal functor is a map $\mathcal{M}^\otimes\to \mathcal{M}^{'\otimes}$ of coCartesian fibrations
over $I$,
which carries coCartesian edges to coCartesian edges.
Let $\widehat{\mathcal{C}\textup{at}}_{\infty}^{\Delta,\textup{sMon}}$
be the simplicial category of symmetric monoidal $\infty$-categories
in which morphisms are symmetric monoidal functors.
Let $\catsmon$ 
be the simplicial nerve of $\widehat{\mathcal{C}\textup{at}}_{\infty}^{\Delta,\textup{sMon}}$ (see \cite[2,1.4.13]{HA}).
For a symmetric monoidal $\infty$-category $\mathcal{M}$, we let $\widehat{\mathcal{C}\textup{at}}_{\infty,\mathcal{M}/}^{\textup{sMon}}$ be
the undercategory in the obvious manner.
We shall refer to a morphism in
$\widehat{\mathcal{C}\textup{at}}_{\infty,\DQ^{\otimes}(k)/}^{\textup{sMon}}$
as a $k$-linear symmetric monoidal functor.
For an affine $k$-scheme $S$, we denote by $\DQ^{\otimes}(S)$
the $\infty$-category $\DQ(S)$ endowed with the natural symmetric
monoidal structure.
If we adopt the notation in Section \ref{subsec:2-2},
the symmetric monoidal $\infty$-category $\DQ^{\otimes}(S)$
is obtained from the pair $(\NNNN(\textup{QC}(S)^c),W_S)$
by applying the left adjoint
functor defined in \cite[preliminary discussion of 4.1.3.4]{HA} to it
(see also \cite[4.1.3.6]{HA}).
Let $\XX$ be an algebraic stack over a field $k$
and let $J$ be the category of affine $k$-schemes over $\XX$.
By this left adjoint functor, the functor given by $\{S\to \XX\} \mapsto (\NNNN(\textup{QC}(S)^c),W_S)$
induces a functor
$J^{op}\to \widehat{\mathcal{C}\textup{at}}_{\infty,\DQ^{\otimes}(k)/}^{\textup{sMon}}$
sending $S \to \XX\in J$ to $\DQ^{\otimes}(S)$
equipped with the structure map $\DQ^\otimes(k)\to \DQ^\otimes(S)$
and sending $f:S'\to S$ to $f^*:\DQ^\otimes(S)\to \DQ^{\otimes}(S')$.
Take a limit $(J^{op})^{\triangleleft}\to \widehat{\mathcal{C}\textup{at}}_{\infty,\DQ^{\otimes}(k)/}^{\textup{sMon}}$ of the diagram
and we shall denote the limit by $\DQ^{\otimes}(\XX)$.
The $\infty$-category $\lcatinf$ endows with a symmetric monoidal
structure given by Cartesian product $\mathcal{C}\times\mathcal{D}$
\cite[2.4]{HA} and a symmetric monoidal $\infty$-category
can be viewed as a commutative algebra (monoid) object in $\lcatinf$.
Thus $\catsmon$ is equivalent to the $\infty$-category
of commutative algebra objects in $\lcatinf$.
By applying \cite[3.2.2.5]{HA} and \cite[1.2.13.8]{HA}
to $\lcatinf$,
the underlying category of $\DQ^{\otimes}(\XX)$ is equivalent to
$\DQ(\XX)$.
If $\mathcal{N},\mathcal{N}'\in \widehat{\mathcal{C}\textup{at}}_{\infty,\mathcal{M}/}^{\textup{sMon}}$
we shall denote by $\Map_{\mathcal{M}}^\otimes(\mathcal{N},\mathcal{N}')$
the mapping space from $\mathcal{N}$ to $\mathcal{N}'$ in
$\widehat{\mathcal{C}\textup{at}}_{\infty,\mathcal{M}/}^{\textup{sMon}}$.
If $\mathcal{M}$ is the initial object, we write
$\Map^\otimes(\mathcal{N},\mathcal{N}')$
for $\Map_{\mathcal{M}}^\otimes(\mathcal{N},\mathcal{N}')$.
We usually write
$\Map_{\DDD(k)}^\otimes(\DQ(\XX), \DQ(S))$
for $\Map_{\DDD^\otimes(k)}^\otimes(\DQ^\otimes(\XX), \DQ^\otimes(S))$.
Let $\Map_{k}^\otimes(\DQ^{\otimes}(\XX),\DQ^{\otimes}(S))$
(or simply $\Map_{k}^\otimes(\DQ(\XX),\DQ(S))$)
denote the full subcategory of
$\Map_{\DDD(k)}^\otimes(\DQ(\XX), \DQ(S))$,
spanned by colimit-preserving functors.
We can regard $\Map_{\DDD(k)}^\otimes(\DQ(\XX),\DQ(S))$ as
the homotopy fiber product of
\[
\Map_{\widehat{\mathcal{C}\textup{at}}_{\infty}^{\Delta,\textup{sMon}}}(\DQ(\XX),\DQ(S))\to \Map_{\widehat{\mathcal{C}\textup{at}}_{\infty}^{\Delta,\textup{sMon}}}(\DQ(k),\DQ(S))\leftarrow \{\ast\},
\]
where
the first map is induced by the structure map $\DQ^\otimes(k)\to \DQ^\otimes(\XX)$
and the second map is induced by the structure map $\DQ^\otimes(k)\to \DQ^\otimes(S)$.
Let $\Hom_k(S,\XX)$ denote a 1-groupoid of $k$-morphisms
from a fixed affine $k$-scheme $S$ to an algebraic stack $\XX$.
We shall regard $\Hom_k(S,\XX)$ as the (simplicial) nerve
of $\Hom_k(S,\XX)$.
Then there is the natural map $\Hom_k(S,\XX)^{op}\to J^{op}$ and
it extends to a map of left cones (cf. \cite{HTT}):
$(\Hom_k(S,\XX)^{op})^{\triangleleft}\to (J^{op})^{\triangleleft}$.
By composing $(J^{op})^{\triangleleft}\longrightarrow \widehat{\mathcal{C}\textup{at}}_{\infty,\DQ^{\otimes}(k)/}^{\textup{sMon}}$ with $(\Hom_k(S,\XX)^{op})^{\triangleleft}\to (J^{op})^{\triangleleft}$ we obtain
\[
\FFF:\Hom_k(S,\XX)\to \Map_{k}^{\otimes}(\DQ^{\otimes}(\XX),\DQ^{\otimes}(S))
\]
which carries $f:S\to \XX$ to $f^*$.

Let $\XX$ be an algebraic stack that satisfies either (i) or (ii)
in Section \ref{subsec:2-3}.
In virtue of Corollary~\ref{summary4} any $k$-linear symmetric
monoidal colimit-preserving functor $\DQ^\otimes(\XX)\to \DQ^\otimes(S)$
preserves vector bundles.
Thus there is a diagram
\[
\xymatrix{
 &  \Map_k^{\otimes}(\DQ^{\otimes}(\XX),\DQ^{\otimes}(S)) \ar[d] \\
 \Hom_{k}(S,\XX)\ar[r]^(0.4){\FFF'} \ar[ru]^{\FFF}& \Map^{\otimes}_{k}(\textup{N}(\textup{h}\DV^{\otimes}(\XX)),\DV^{\otimes}(S))\\
}
\]
in $\lcatinf$, where the vertical arrow is induced by the restriction,
 $\FFF'$ sends $f:S\to \XX$ to the $k$-linear symmetric monoidal
functor
$f^*:\textup{N}(\textup{h}\DV^{\otimes}(\XX))\to\DV^{\otimes}(S)$, and
$\Map^{\otimes}_{k}(\textup{N}(\textup{h}\DV^{\otimes}(\XX)),\DV^{\otimes}(S))$
is the category of $k$-linear symmetric monoidal additive exact functors
in which morphisms are monoidal natural transformations.

The following is Tannaka duality for quasi-projective schemes
with action of an affine group scheme (generalizing
the classical one), proved by Savin \cite{Sav}.
In \cite{LT}, Lurie shows another version of Tannaka duality for a geometric stack using
the symmetric monoidal category of quasi-coherent sheaves.

\begin{Theorem}[\cite{Sav}]
\label{full}
Suppose that $\XX$ is a quotient stack of the form $[X/G]$,
where $X$ is a separated noetherian scheme and $G$ is a linear algebraic group
acting on $X$.
Suppose that there is a very ample $G$-invertible sheaf on $X$.
The functor $\FFF'$ is a categorical equivalence.
\end{Theorem}

\vspace{2mm}

Next we generalize extension lemmas proved in Section \ref{sec:3}
to a version of symmetric monoidal functors (Proposition~\ref{monoidalextension}).
We shall begin by describing the naive idea of this generalization.
Let $\mathcal{A}^\otimes$ and $\mathcal{B}^\otimes$
be two symmetric monoidal $\infty$-categories.
Let $\mathcal{A}_c^\otimes$ be a
symmetric monoidal full subcategory of $\mathcal{A}^\otimes$.
(To simplify the problem, we may suppose further that $\infty$-categories
$\mathcal{A}$
and $\mathcal{B}$ are $1$-categories.) 
Note that a symmetric monoidal $\infty$-category $\mathcal{A}^\otimes$
amounts to
$f_{\mathcal{A}^\otimes}:I\to \widehat{\mathcal{C}\textup{at}}_\infty$
such that
$\alpha^{j,i}_{!}:f(\langle i\rangle_\ast)\to f(\langle 1\rangle_\ast)$
($1\le j\le i$) induces an equivalence $f(\langle i\rangle_\ast)\to f(\langle 1\rangle_\ast)\times\ldots\times f(\langle 1\rangle_\ast)$ to the $i$-fold product,
that is, a commutative monoid object.
More precisely, $\catsmon$ can be embedded into $\Fun(I,\widehat{\mathcal{C}\textup{at}}_{\infty})$ as the full subcategory spanned by commutative monoid
objects (see \cite[2.4.2.6]{HA}).
Informally, $f_{\mathcal{A}^\otimes}:I\to \lcatinf$ is depicted as
\[
\xymatrix@R=5mm @C=0mm{
\cdots & \rightleftharpoons & \mathcal{A}^{\times i}& \rightleftharpoons & \mathcal{A}^{\times i-1} & \rightleftharpoons & \cdots & \rightleftharpoons & \mathcal{A}^{\times 2} & \rightleftharpoons &  \mathcal{A} & \rightleftharpoons &\mathcal{A}^{\times 0}\simeq\Delta^0  \\
}
\]
where
$\rightleftharpoons$ between $\mathcal{A}^{\times i}$ and $\mathcal{A}^{\times i-1}$ informally
represents morphisms induced by maps between $\langle i\rangle_*$ and $\langle i-1\rangle_*$ (namely, $\mathcal{A}^{\times i}$ is $f_{\mathcal{A}^\otimes}(\langle i\rangle_*)$).
A symmetric monoidal functor amounts to a natural transformation
$I \times \Delta^1\to \widehat{\mathcal{C}\textup{at}}_{\infty}$
between commutative monoid objects.
It is informally described by the diagram in $\lcatinf$:
\[
\xymatrix@R=5mm @C=0mm{
\cdots & \rightleftharpoons &\mathcal{A}^{\times i} \ar[d] & \rightleftharpoons & \mathcal{A}^{\times i-1} \ar[d] & \rightleftharpoons &\cdots & \rightleftharpoons &\mathcal{A}^{\times 2} \ar[d] & \rightleftharpoons &  \mathcal{A} \ar[d] & \rightleftharpoons &\mathcal{A}^{\times0} \ar[d] \\
\cdots & \rightleftharpoons &\mathcal{B}^{\times i}& \rightleftharpoons & \mathcal{B}^{\times i-1} & \rightleftharpoons &\cdots & \rightleftharpoons &\mathcal{B}^{\times 2} & \rightleftharpoons &  \mathcal{B} & \rightleftharpoons &\mathcal{B}^{\times0}.
}
\]
Let $\Map_{\star}^\otimes(\mathcal{A}^\otimes,\mathcal{B}^\otimes)$
be a full subcategory, and suppose that we want to prove that $\Map_{\star}^\otimes(\mathcal{A}^\otimes,\mathcal{B}^\otimes) \to \Map^\otimes(\mathcal{A}_c^\otimes,\mathcal{B}^\otimes)$ induced by the inclusion $\mathcal{A}_c^\otimes \subset \mathcal{A}^\otimes$ is fully faithful.
The rough idea is to concentrate on a full subcategory of 
$\widehat{\mathcal{C}\textup{at}}_\infty$ consisting of objects
of the image of $f_{\mathcal{A}^\otimes}$, $f_{\mathcal{A}_c^\otimes}$
and $f_{\mathcal{B}^\otimes}$ and restrict
morphisms from $\mathcal{A}_c^{\times i}$ to those morphisms
which can be uniquely extended
to morphisms from $\mathcal{A}^{\times i}$ in an appropriate sense (involving
the subscript ``$\star$'').
To explain this,
let us assume that we have a subcategory $H$
of $\widehat{\mathcal{C}\textup{at}}_{\infty}$ having the following properties:
\begin{itemize}
\item objects of $H$ are $\{\Delta^0, \mathcal{A},\mathcal{A}^{\times 2},\ldots,\mathcal{A}^{\times i},\ldots\}\cup
\{\Delta^0, \mathcal{B},\mathcal{B}^{\times 2},\ldots,\mathcal{B}^{\times i},\ldots\}$,
\item any morphism $\mathcal{A}^{\times i}\to \mathcal{A}^{\times j}$ in $H$
carries $\mathcal{A}_c^{\times i}$ to $\mathcal{A}_c^{\times j}$,

\item if we denote by $\Map_H(-,-)$ the mapping space in $H$,
$\Map_H(\mathcal{A}^{\times i},\mathcal{A}^{\times j}) \to \Map(\mathcal{A}_c^{\times i},\mathcal{A}_c^{\times j})$ is fully faithful for any $i,j\ge 0$,

\item $\Map_H(\mathcal{A}^{\times i},\mathcal{B}^{\times j})\to \Map(\mathcal{A}^{\times i}_c,\mathcal{B}^{\times j})$ is fully faithful for any $i,j\ge 0$,

\item $\Map_H(\mathcal{B}^{\times i},\mathcal{A}^{\times j})$ is the
empty set for any $i\ge 0$ and $j\ge 1$,

\item $f_{\mathcal{A}^\otimes},f_{\mathcal{B}^\otimes}:I\rightrightarrows \widehat{\mathcal{C}\textup{at}}_{\infty}$ can factor through $H$,

\item if $I\times \Delta^1 \rightarrow \widehat{\mathcal{C}\textup{at}}_{\infty}$ corresponds to an object in $\Map_{\star}^\otimes(\mathcal{A}^\otimes,\mathcal{B}^\otimes)$, it factors through $H$.

\end{itemize}
Let $H'$ be the full subcategory of $\widehat{\mathcal{C}\textup{at}}_{\infty}$ consisting of
\[
\{\Delta^0, \mathcal{A}_c,\mathcal{A}_c^{\times 2},\ldots,\mathcal{A}_c^{\times i},\ldots\}\cup \{\Delta^0, \mathcal{B},\mathcal{B}^{\times 2},\ldots,\mathcal{B}^{\times i},\ldots\}.
\]
Let $z:H\to H'$ be the functor determined by the restrictions $\mathcal{A}_c^{\times i}\subset \mathcal{A}^{\times i}$,
which carries $\mathcal{A}^{\times i}$
and $\mathcal{B}^{\times i}$ to
$\mathcal{A}_c^{\times i}$
and $\mathcal{B}^{\times i}$ respectively.
Then by the above properties, for any two $X,Y\in H$, $z$ induces a fully faithful functor
$\Map_H(X,Y)\to \Map(z(X),z(Y))$.
Observe that this faithfulness implies that 
$\Map_{\star}^\otimes(\mathcal{A}^\otimes,\mathcal{B}^\otimes) \to \Map^\otimes(\mathcal{A}_c^\otimes,\mathcal{B}^\otimes)$ is fully faithful.

We apply this idea to prove Proposition~\ref{monoidalextension}.
For this purpose, we will define two simplicial categories $\mathscr{P}$
and $\mathscr{Q}$.
Here we use notation similar to Section \ref{sec:3}, i.e.,
$\CCC=\DQ(\XX),\ \DDD=\DQ(S),\ \KKK=\DQ(k)$.
(We here need to take the $k$-linear structures into consideration
and consider three symmetric monoidal $\infty$-categories.)
Let $\alpha,\beta,\gamma:I\to \widehat{\mathcal{C}\textup{at}}_{\infty}$
be functors corresponding to $\CCC^\otimes$, $\DDD^\otimes$ and
$\KKK^\otimes$ respectively.
A symmetric monoidal functor $\DQ^{\otimes}(\XX)\to \DQ^{\otimes}(S)$ amounts
to a natural transformation $\alpha\to \beta$, that is, a morphism from $\alpha$ to $\beta$ in $\Fun (I,\lcatinf)$.

We define a fibrant simplicial category $\mathscr{C}$.
Objects of $\mathscr{C}$ are $\CCC^{\times i},\DDD^{\times i}$ and $\KKK^{\times i}$ ($i\ge0 $).
For simplicity, $\KKK^{\times i}=\mathcal{A}_0^i$,
$\CCC^{\times i}=\mathcal{A}_1^i$ and
$\DDD^{\times i}=\mathcal{A}_2^i$.
For $0\le r,s\le 2$,
a simplicial set $\Map_{\mathscr{C}}(\mathcal{A}_r^i,\mathcal{A}_s^j)$
is $\Map(\mathcal{A}_r^i,\mathcal{A}_s^j)$.

Let $\mathcal{V}^i_r$ denote the full subcategory of $\mathcal{A}_r^i$ which
consists of vector bundles, that is, objects belonging to
$((\mathcal{A}_r)_{\textup{vect}})^{\times i}$.
Here $(\mathcal{A}_r)_{\textup{vect}}$ denotes
 the full subcategory of $\mathcal{A}_r$ spanned by complexes which are
 quasi-isomorphic to vector bundles.
Let $K$ be a small simplicial set.
Let $K\to \mathcal{A}^{i}_r$ be a functor which has the property:
there is $c\in \{1,\ldots,i\}$ such that for $j\neq c$
the composite $K\to \mathcal{A}_r^1$
with the $j$-th projection
$\alpha^{j,i}_{!}:\mathcal{A}_r^i\to \mathcal{A}_r^1$ is equivalent to
a constant diagram. We will call such a functor
a diagram of one-variable indexed by $K$.
We say that a functor $\mathcal{A}_r^i\to \mathcal{A}_s^j$
is good if for any simplicial set $K$
and any diagram of one-variable $K\to \mathcal{A}_r^i$
the composite $K\to \mathcal{A}_s^j$ is a diagram of one-variable.

Let us
define a simplicial subcategory $\mathscr{P}$ of $\mathscr{C}$ as follows.
A collection of objects of $\mathscr{P}$ is the same as that of $\mathscr{C}$.
We define hom simplicial sets as follows:
For any $0\le r,s\le 2$,
$\Map_{\mathscr{P}}(\mathcal{A}_r^i,\mathcal{A}_s^j)$
is the full subcategory of
$\Map(\mathcal{A}_r^i,\mathcal{A}_s^j)$,
spanned by functors satisfying the properties:
\begin{itemize}
\item good,

\item colimit-preserving separately in each variable
of $\mathcal{A}_r^i$,

\item it sends $\mathcal{V}_r^j$ to $\mathcal{V}_s^j$.

\end{itemize}
These data constitute a simplicial category $\mathscr{P}$.
Next we will define another fibrant simplicial category $\mathscr{Q}$.
The collection of objects is the same as $\mathscr{P}$.
For any $0\le r,s\le 2$,
$\Map_{\mathscr{Q}}(\mathcal{A}_r^i,\mathcal{A}_s^j)$
is the full subcategory of
$\Map(\mathcal{V}_r^i,\mathcal{V}_s^j)$,
spanned by functors which are
equivalent to image of the restriction
$\Map_{\mathscr{P}}(\mathcal{A}_r^i,\mathcal{A}_s^j)\to \Map(\mathcal{V}_r^i,\mathcal{V}_s^j)$. Compositions are well-defined and the data
form a simplicial category.
There is a natural simplicial functor $\xi:\mathscr{P}\to \mathscr{Q}$.
Then by Lemma~\ref{extensionstep1} and~\ref{extensionstep5},
we have the following:

\begin{Lemma}
Let $\XX$ and $S$ be perfect stacks over $k$.
Suppose that $\XX$ and $S$ have cohomological dimension zero.
 Then the simplicial functor $\xi:\mathscr{P}\to \mathscr{Q}$
 of fibrant simplicial categories is an equivalence.\end{Lemma}

Since $\alpha,\beta$ and $\gamma$ factor through $\textup{N}(\mathscr{P})$,
we write $\alpha',\beta',\gamma':I\to \textup{N}(\mathscr{P})$ for their factorizations.
A symmetric monoidal functor $\DQ^\otimes (\XX)\to \DQ^\otimes (S)$
can be viewed as a natural transformation
$\alpha'\to \beta'$.
On the other hand, a natural transformation between
$\textup{N}(\xi)\circ\alpha',\textup{N}(\xi)\circ\beta':I\to \textup{N}(\mathscr{P})\to \textup{N}(\mathscr{Q})$
is nothing but a symmetric monoidal functor between
$\DV^\otimes(\XX)$ and $\DV^\otimes(S)$.
By these observations we deduce the following:

\begin{Proposition}
\label{monoidalextension}
Suppose that $\XX$ has cohomological dimension zero $($and $S$
is affine$)$.
Let $\Phi,\Psi:\DQ^{\otimes}(\XX)\to \DQ^{\otimes}(S)$ be $k$-linear
symmetric monoidal
functors
which preserve colimits.
Let $\overline{\Phi}:\DV^{\otimes}(\XX)\to \DV^{\otimes}(S)$
 and $\overline{\Psi}:\DV^{\otimes}(\XX)\to \DV^{\otimes}(S)$ be the restriction of
$\Phi$ and $\Psi$. $($By Corollary~\ref{summary4},
$\Phi$ and $\Psi$ preserve vector bundles.$)$
Then the natural functor
\[
\Map(\Phi,\Psi)\to \Map(\overline{\Phi},\overline{\Psi})
\]
is a weak homotopy equivalence.
Here we denote by $\Map(\Phi,\Psi)$ and $\Map(\overline{\Phi},\overline{\Psi})$
the mapping spaces in
$\Map_k^\otimes(\DQ(\XX),\DQ(S))$
and $\Map_{\DV^{\otimes}(k)}^{\otimes}(\DV^\otimes(\XX),\DV^\otimes(S))$
respectively.
\end{Proposition}

\begin{Remark}
Although it is sufficient for our main goal,
we impose the unpleasant condition on cohomological
dimension in Proposition~\ref{monoidalextension}
as well as Proposition~\ref{cofinality}.
It is desirable to remove this condition;
it is  meaningful for other applications to generalize
Proposition~\ref{cofinality} (and Proposition~\ref{monoidalextension})
to the case of noetherian stacks.
\end{Remark}

\begin{Remark}
\label{monoidalqcperf}
Using an argument which is similar to Proposition~\ref{monoidalextension}
and Lemma~\ref{extensionstep1},
we deduce that there is an natural equivalence
\[
\Map^{\otimes}_k(\DQ(\XX),\DQ(S))\to 
\Map^{\otimes}_k(\DP(\XX),\DP(S)).
\]
where $\Map^{\otimes}_k(\DP(\XX),\DP(S))$
denotes the mapping space of symmetric monoidal exact functors
$\DP^\otimes(\XX)\to \DP^\otimes(S)$ over $\DP^\otimes(k)$.
(Note that any symmetric monoidal functor preserves dualizable objects,
i.e., perfect complexes.)
Therefore we can replace $\Map^{\otimes}_k(\DQ^\otimes(\XX),\DQ^\otimes(S))$
by $\Map^{\otimes}_k(\DP^\otimes(\XX),\DP^\otimes(S))$
in Theorem~\ref{supermain}.
\end{Remark}

\begin{Lemma}
\label{symmetricopenreduction2}
The restriction functor $\DQ^{\otimes}(\XX)\to \DQ^{\otimes}(\mathcal{U})$
induces a categorical equivalence
\[
\Map^{\otimes}(\DQ(\mathcal{U}),\DQ(S))\to
\Map^{\otimes}_{\diamond}(\DQ(\XX),\DQ(S)).
\]
The notation $\Map^{\otimes}_{\diamond}(\DQ(\XX),\DQ(S))$
indicates the full subcategory spanned by functors $\Phi$
such that
if $f:H\to G\in \Fun(\Delta^1,\DQ(\XX))$ induces
an equivalence in $\DQ(\UU)$ then $\Phi(f)$ is an equivalence.
\end{Lemma}

\Proof
Let $p:\CCC^\otimes\to I$ and 
$q:\DDD^\otimes \to I$
be coCartesian fibrations
that correspond to the symmetric monoidal $\infty$-categories
$\DQ^\otimes(\XX)$ and $\DQ^\otimes (S)$.
Let $\UU\subset \XX$ be an open substack.
Let $\CCC_{\UU}^\otimes\subset \CCC^\otimes$
be the full subcategory 
such that $\CCC_{\UU}^\otimes\cap p^{-1}(\langle n\rangle_*)$
is spanned by $\DQ(\UU)\times\cdots \times\DQ(\UU)$ ($n$-times product).
Namely, $p_{\UU}:\CCC_{\UU}^\otimes\to I$ is a coCartesian fibration
that corresponds to $\DQ^\otimes(\UU)$.
Let $\Phi:\CCC^\otimes \to \DDD^\otimes$ be a symmetric monoidal
functor
such that
if $f:H\to G\in \Fun(\Delta^1,\DQ(\XX))$ induces
an equivalence in $\DQ(\UU)$ then $\Phi(f)$ is an equivalence.
To prove our claim, in the light of \cite[4.3.2.15]{HTT} and Lemma~\ref{openreduction}
it will suffice to show that
$\Phi(P)$ is
a $q$-limit of
the diagram $(\CCC_{\UU}^\otimes)_{P/}\to \CCC^\otimes\to \DDD^\otimes$
for any $P\in \CCC^\otimes$ where $(\CCC_{\UU}^\otimes)_{P/}$
denotes the undercategory.
Suppose that $P\in p^{-1}(\langle n\rangle_*)$ and 
$P=[P_1,\ldots,P_n]\in \DQ(\XX)^{\times n}$.
Let $P_{\UU}$ be
$[L(P_1),\ldots,L(P_n)]\in \CCC_{\UU}^\otimes$ where
$L:\DQ(\XX)\to \DQ(\UU)\to \DQ(\XX)$ (cf. Lemma~\ref{openreduction}).
Here we refer to $P_{\UU}$ as a $\UU$-localization of $P$.
Since a $\UU$-localization of $\iota_{!}(P_{\UU})$ is equivalent
to a $\UU$-localization of $\iota_!(P)$
for any $\iota\in \Fun(\Delta^1,I)$ and $p$ is a coCartesian fibration,
we see that $P\to P_{\UU}$ is an initial object
of $\CCC_{\UU}^\otimes\times_{\CCC^\otimes}\CCC^{\otimes}_{P/}$
(cf. \cite[5.2.7.6]{HTT}).
Thus unwinding the definition of $q$-limits \cite[4.3.1.1]{HTT}
we conclude that $\Phi(P)$ is a $q$-limit of
$(\CCC_{\UU}^\otimes)_{P/}\to \CCC^\otimes\to \DDD^\otimes$.
\QED

Let $\DQ^\otimes(k)\to \DQ^\otimes(\UU)$ and $\DQ^\otimes(k)\to \DQ^\otimes(\XX)$
be $k$-linear structure maps.
These maps induce
\[
u:\Map^{\otimes}(\DQ(\mathcal{U}),\DQ(S))\to \Map^{\otimes}(\DQ(k),\DQ(S))
\]
and
\[v:\Map^{\otimes}_{\diamond}(\DQ(\mathcal{X}),\DQ(S))\to \Map^{\otimes}(\DQ(k),\DQ(S)).
\]
Let $\iota:\Delta^0=\ast\to \Map^{\otimes}(\DQ(k),\DQ(S))$ be
the map corresponds to $\DQ^\otimes (k)\to \DQ^\otimes (S)$.
Note that the mapping space $\Map^{\otimes}_{\DDD(k)}(\DQ(\mathcal{U}),\DQ(S))$
is a homotopy pullback of $u$ along $\iota$.
We shall denote by $\Map^{\otimes}_{\DDD(k),\diamond}(\DQ(\XX),\DQ(S))$ a homotopy pullback of $v$ along $\iota$.

\begin{Corollary}
There is a natural equivalence
\[
\Map^{\otimes}_{\DDD(k)}(\DQ(\mathcal{U}),\DQ(S))\to
\Map^{\otimes}_{\DDD(k),\diamond}(\DQ(\XX),\DQ(S)).
\]
\end{Corollary}

We would like to record a direct consequence of Lemma~\ref{ff;1cat}.

\begin{Lemma}
\label{symmetrictruncation}
The natural equivalence $\DV^{\otimes}(\XX)\stackrel{\sim}{\to} \textup{N}(\textup{h}\DV^{\otimes}(\mathcal{X}))$
induces a categorical equivalence
\[
\Map^{\otimes}(\textup{N}(\textup{h}\DV(\mathcal{X})),\DV(S))\to
\Map^{\otimes}(\DV(\XX),\DV(S)).
\]
Moreover the equivalence induces an equivalence
\[
\Map^{\otimes}_{\DV^\otimes(k)}(\textup{N}(\textup{h}\DV(\mathcal{X})),\DV(S))\to
\Map^{\otimes}_{\DV^\otimes(k)}(\DV(\XX),\DV(S)).
\]
\end{Lemma}

The latter follows from the fact:
$\Map^{\otimes}_{\DV^\otimes(k)}(\textup{N}(\textup{h}\DV(\mathcal{X})),\DV(S))$ is a homotopy pullback of the map
$\Map^{\otimes}(\textup{N}(\textup{h}\DV(\mathcal{X})),\DV(S))
\to \Map^\otimes(\DV(k),\DV(S))$ (induced by $\DV^\otimes(k)\to
\textup{N}(\textup{h}\DV(\mathcal{X}))$) along $\ast\to
\Map^\otimes(\DV(k),\DV(S))$ which is induced by the $k$-linear structure map
$\DV^\otimes(k)\to \DV^\otimes(S)$.
In a similar way, 
$\Map^{\otimes}_{\DV^\otimes(k)}(\DV(\mathcal{X}),\DV(S))$ is a homotopy pullback of
$\Map^{\otimes}(\DV(\mathcal{X}),\DV(S))
\to \Map^\otimes(\DV(k),\DV(S))$
along $\ast\to
\Map^\otimes(\DV(k),\DV(S))$.

\begin{Proposition}
\label{combining}
Let $\XX$ be an algebraic stack over $k$
that satisfies either (i) or (ii)
in Section \ref{subsec:2-3}.
Let $S$ be an affine scheme over $k$.
Let $a,b:\DQ^\otimes(\XX)\to \DQ^{\otimes}(S)$
be $k$-linear symmetric monoidal colimit-preserving functors
$($such that $a(\DV(\XX))$ and $b(\DV(\XX))$ lie in $\DV(S)$$)$.
Let $\bar{a}:=a|_{\DV^\otimes(\XX)}:\DV^\otimes(\XX)\to \DV^\otimes(S)$ and
$\bar{b}:=b|_{\DV^\otimes(\XX)}:\DV^\otimes(\XX)\to \DV^\otimes(S)$. 
Let $\Map^{\otimes}(a,b)$ be the mapping space from $a$ to $b$ in
$\Map^{\otimes}_{k}(\DQ(\XX),\DQ(S))$
and let $\Map(\bar{a},\bar{b})$ be the mapping space from $\bar{a}$
to $\bar{b}$ in
$\Map_{\DV^{\otimes}(k)}^\otimes(\DV(\XX),\DV(S))$.
Suppose that there exist $x,y:S\rightrightarrows \XX$
such that $\bar{a}$ and $\bar{b}$ are equivalent to pullback
functors $x^*:\DV(\XX)\to \DV(S)$ and $y^*:\DV(\XX)\to \DV(S)$ as 
functors respectively.
Then the restriction $\Map^{\otimes}(a,b)\to \Map^{\otimes}(\bar{a},\bar{b})$
is a weak homotopy equivalence.
\end{Proposition}

\Proof
We first fix some notation.
Take a Zariski affine covering $\sqcup_l S_l\to S$
such that each $S_l\to S\stackrel{x}{\to} \XX$
(and $S_l\to S\stackrel{y}{\to} \XX$)
factors through a quasi-compact open substack
$\UU_l\subset \XX$ which has cohomological dimension zero.
Let $Z$ be the category of affine schemes $T$ over $S$
such that $T\to S$ is an open immersion
and $T\to S$ factors through some $S_l\subset S$.
Then $\DQ^\otimes(S)$ is a limit of
the diagram $Z^{op}\to \widehat{\mathcal{C}\textup{at}}_{\infty,\DQ^\otimes(k)/}^{\textup{sMon}}$ sending $T\in Z$ to $\DQ^\otimes(T)$.
If $z\in Z$ indicates $T\to S$,
then we write $S^z$ for $T$ and we denote by $a^z$ the composite
$\DQ^\otimes(\XX)\stackrel{a}{\to} \DQ^\otimes(S)\to \DQ^\otimes(S^z)$, and we use the notation $b^z$
and $\bar{a}^z$ in a similar manner.
Note that by Lemma~\ref{bundlecolimit} and Remark~\ref{sketo}
a bounded complex $P$ on $\XX$ which in each degree is an infinite
direct sum of vector bundles is a finite colimit of
infinite direct sums of vector bundles up to shifts.
Thus $a(P)\simeq x^*(P)$ and $b(P)\simeq y^*(P)$.
If $P\in \DQ(\XX)$ is acyclic on $\UU_l$,
then by Lemma~\ref{inductive} we conclude that $a(P)$ and $b(P)$
are acyclic on $S_l$.
Thus by Lemma~\ref{symmetricopenreduction2}
each $\DQ^\otimes(\XX)\stackrel{a}{\to} \DQ^\otimes(S)\to\DQ^\otimes(S_l)$ and $\DQ^\otimes(\XX)\stackrel{b}{\to} \DQ^\otimes(S)\to \DQ^\otimes(S_l)$
factor through $\DQ^\otimes(\XX)\to \DQ^\otimes(\UU_l)$.
To prove our claim, it is convenient to recall the presentation of
the mapping spaces in an $\infty$-category, introduced in \cite[page 28]{HTT}.
Let $\mathcal{C}$ be a Kan complex (for the case of $\infty$-categories see
\cite{HTT}) and let $c$ and $c'$
be two objects in $\CCC$, i.e., two vertices.
We define a mapping Kan complex $\Fun^{(c,c')}(\Delta^1,\mathcal{C})$
to be the fiber product of $\Fun(\Delta^{1},\CCC)\to \Fun(\partial\Delta^1,\CCC)\leftarrow \ast=\{(c,c')\}$. Since $\Fun(\Delta^{1},\CCC)\to \Fun(\partial\Delta^1,\CCC)$ is a Kan fibration induced by inclusion
$\partial\Delta^1\to \Delta^1$, the fiber product is a homotopy fiber product.
Using this presentation and the universality of limits together with
Lemma~\ref{symmetricopenreduction2} and Remark~\ref{monoidalqcperf},
we have the following categorical equivalences
\begin{equation*}
\begin{split}
\Fun^{(a,b)}(\Delta^1,\Map^{\otimes}_{\DDD(k)}(\DQ(\XX),&\DQ(S))) \\
&\simeq \Fun^{(a,b)}(\Delta^1,\lim_{z\in Z}\Map^{\otimes}_{\DDD(k)}(\DQ(\XX),\DQ(S^z))) \\
& \simeq \lim_{z\in Z}\Fun^{(a^z,b^z)}(\Delta^1,\Map^{\otimes}_{\DDD(k)}(\DQ(\XX),\DQ(S^z)) \\
& \simeq \lim_{z\in Z}\Fun^{(a^z,b^z)}(\Delta^1,\Map^{\otimes}_{\DDD(k)}(\DQ(\sqcup_l\ \UU_l),\DQ(S^z))) \\
& \simeq \lim_{z\in Z}\Fun^{((a')^z,(b')^z)}(\Delta^1,\Map^{\otimes}_{k}(\DP(\sqcup_l\ \UU_l),\DP(S^z)).
\end{split}
\end{equation*}
Here $(a')^z$ and $(b')^z$ is the restriction of $a^z$ and $b^z$
to $\DP^\otimes (\XX)$ respectively.
Applying Proposition~\ref{monoidalextension} we have
\begin{equation*}
\begin{split}
\lim_{z\in Z}\Fun^{((a')^z,(b')^z)}(\Delta^1,\Map^{\otimes}_{k}&(\DP(\sqcup_l\ \UU_l),\DP(S^z)))\\
&\simeq \lim_{z\in Z}\Fun^{(\bar{a}^z,\bar{b}^z)}(\Delta^1,\Map^{\otimes}_{\DV^{\otimes}(k)}(\DV(\sqcup_l\ \UU_l),\DV(S^z))) \\
&\simeq \lim_{z\in Z}\Fun^{(\bar{a}^z,\bar{b}^z)}(\Delta^1,\Map^{\otimes}_{\DV^\otimes(k)}(\textup{h}\DV(\sqcup_l\ \UU_l),\DV(S^z))).
\end{split}
\end{equation*}
We abusively write $\textup{h}\DV(\bullet)$ for
$\textup{N}(\textup{h}\DV(\bullet))$.
Note that by Theorem~\ref{full}
for a Zariski open substack $\mathcal{U}\subset \XX$
the full subcategory of
$\Map^{\otimes}_{\DV^{\otimes}(k)}(\textup{h}\DV(\UU),\DV(S^z))$,
spanned by additive exact functors
can be viewed as $\Hom_k(S^z,\UU)$. Thus
the full subcategory of
$\Map^{\otimes}_{\DV^{\otimes}(k)}(\textup{h}\DV(\UU),\DV(S^z))$ can be naturally viewed as a
full subcategory of
$\Map^{\otimes}_{\DV^{\otimes}(k)}(\textup{h}\DV(\XX),\DV(S^z))$.
Thus by these obsevations, the descent theory of
vector bundles and Lemma~\ref{symmetrictruncation}
we have equivalences
\begin{equation*}
\begin{split}
\lim_{z\in Z}\Fun^{(\bar{a}^z,\bar{b}^z)}(\Delta^1,\Map^{\otimes}_{\DV^{\otimes}(k)}(\textup{h}&\DV(\sqcup_l\ \UU_l),\DV(S^z)))\\
&\simeq \lim_{z\in Z}\Fun^{(\bar{a}^z,\bar{b}^z)}(\Delta^1,\Map^{\otimes}_{\DV^{\otimes}(k)}(\textup{h}\DV(\XX),\DV(S^z))) \\
&\simeq \Fun^{(\bar{a},\bar{b})}(\Delta^1,\Map^{\otimes}_{\DV^{\otimes}(k)}(\textup{h}\DV(\XX),\DV(S))) \\
&\simeq \Fun^{(\bar{a},\bar{b})}(\Delta^1,\Map^{\otimes}_{\DV^{\otimes}(k)}(\DV(\XX),\DV(S))). 
\end{split}
\end{equation*}
Therefore we obtain the desired equivalence.
\QED

Finally, we obtain our main goal:

\begin{Theorem}
\label{supermain}
Let $\XX$ be an algebraic stack which satisfies either (i) or (ii)
in Section \ref{subsec:2-3}.
Let $S$ be a scheme over $k$
 $($we always assume that $S$ is
quasi-compact and has affine diagonal$)$.
Then there is a categorical equivalence
\[
\FFF:\Hom_k(S,\XX)\longrightarrow \Map_k^\otimes(\DQ^{\otimes}(\XX),\DQ^{\otimes}(S))
\]
which sends $f:S\to \XX$ to $f^*$.
\end{Theorem}

\Proof
If $S$ is affine, our claim follows from Corollary~\ref{summary4}, Proposition~\ref{combining}
and Theorem~\ref{full}.
If $S$ is a scheme, take a Zariski covering $\sqcup T\to S$
by affine schemes. It gives rise to a simplicial scheme $S_\bullet\to S$.
Then $\Hom_k(S,\XX)$ is a limit of the cosimplicial diagram of
$\Hom_k(S_i,\XX)$ indexed by $i\in\Delta$.
On the other hand, $\DQ^\otimes(S)$ is a limit of the cosimplicial
diagram $\DQ^\otimes(S_i)$ in $\widehat{\mathcal{C}\textup{at}}_{\infty, \DQ^{\otimes}(k)/}^{\textup{sMon}}$.
Since $\DQ(\XX)\to\DQ(S)$ preserves small colimits if and only if
the composite $\DQ(\XX)\to \DQ(\sqcup T)$ preserves small colimits,
thus $\Map_k^\otimes(\DQ^{\otimes}(\XX),\DQ^{\otimes}(S))$
is a limit of the cosimplicial diagram of
$\Map_k^\otimes(\DQ^{\otimes}(\XX),\DQ^{\otimes}(S_i))$.
Now our assertion follows from the case where $S$ is affine.
\QED

\begin{Remark}
We would like to explain the reason
why we should employ the theory of $(\infty,1)$-categories.
Note that morphisms to $\XX$ have the descent property.
Namely, if $p:S'\to S$ is an \'etale surjective morphism
and $\textup{pr}_1,\textup{pr}_2:S'\times_SS'\rightrightarrows S'$
are the first and second projections,
then a morphism $f':S'\to \XX$ such that $\textup{pr}_1\circ f'=\textup{pr}_2\circ f'$ descents to a unique morphism $f:S\to \XX$ such that $p\circ f=f'$.
Now suppose that Tannaka
duality formulated with the triangulated categories holds.
Then the descent property of morphisms to $\XX$ implies
that functors $\textup{D}(\XX)\to \textup{D}(S)$ of triangulated 
categories of a certain type have
the descent property, where $\textup{D}(\bullet)$ denotes 
the triangulated category of quasi-coherent complexes (or perfect complexes).
However, we can not hope that
the derived categories have a reasonable descent theory.
One of sources of this problem comes from the fact that triangulated categories forget the structure of homotopy coherence
which naturally arise from (co)chain complexes.
Inspired by the derived algebraic geometry \cite{TV0}, \cite{TV}, \cite{DAG}
and derived Morita theory \cite{To}, \cite{BFN},
in order to establish our Tannaka duality
we use not triangulated categories
but ``enhanced higher categories'' such as stable (symmetric
monoidal) $\infty$-categories.
The idea of usage of higher category theory could be found in
algebraic K-theory \cite{TVK}.

We here call Theorem~\ref{supermain} {\it derived Tannaka duality},
which is a title of this section.
But perhaps it is more appropriate to say that Theorem~\ref{supermain}
 is a {\it stable} analogue of Tannaka duality, although the term ``stable analogue'' is ambiguous as well
as the term ``derived analogue''.
\end{Remark}

Let us consider the ($\infty$-)stack
on the \'etale site $(\textup{Aff}_k)$ of affine $k$-schemes:
\[
\mathcal{F}_{\XX}:(\textup{Aff}_k)^{op}\longrightarrow \mathcal{S}
\]
which sends $S$ to $\Map_{k}^{\otimes}(\DQ(\XX),\DQ(S))$.
Here $\mathcal{S}$ is $\infty$-category of spaces (Kan complexes)
\cite[1.2.16]{HTT} and we view $\mathcal{F}_{\XX}$ as an object in $\Fun((\textup{Aff}_k)^{op},
\mathcal{S})$ or the localization of
$\Fun((\textup{Aff}_k)^{op},
\mathcal{S})$ with respect to the \'etale topology of $(\textup{Aff}_k)$
\cite[6.2.2]{HTT}.
An immediate consequence of Theorem~\ref{supermain} is:

\begin{Corollary}
\label{reconstruction}
Let $\XX$ be an algebraic stack over $k$ that satisfies the condition
either (i) or (ii).
Then the stack $\XX$ over $(\textup{Aff}_k)$ is
equivalent to $\mathcal{F}_{\XX}$.
\end{Corollary}

\begin{Remark}
The above corollary is a reconstruction result.
Our reconstruction is of different nature from
one in \cite{Bal}.
The point is that
(i) in loc. cit., schemes are reconstructed as
ringed spaces, whereas we reconstruct them as sheaves on $(\textup{Aff}_k)$
(so it is applicable to the case of stacks),
(ii) on one hand we recover a scheme $X$
from a symmetric monoidal $\infty$-category $\DQ^\otimes(X)$ or $\DP^{\otimes}(X)$;
on the other hand, in loc. cit.,  a scheme is recovered from
a symmetric monoidal triangulated category $\textup{D}_{\textup{perf}}^\otimes(X)$.
We expect that an enhancement of a symmetric monoidal triangulated category $\textup{D}_{\textup{qcoh}}(X)$ is unique in an appropriate sense. In this direction, in the recent paper
\cite{OL} by Lunts and Orlov
it is shown that for a quasi-projective variety $X$
an dg-enhancement of a triangulated category
 $\textup{D}_{\textup{qcoh}}(X)$ is unique.
\end{Remark}

\begin{Remark}
Our result is  also closely related to the moduli of perfect complexes.
Let $X$ be a smooth projective variety.
Let $\mathcal{M}_{\textup{perf}}(X):(\textup{Aff}_k)^{op}\to \mathcal{S}$
be the functor (moreover $\infty$-stack)
which to any $S\in \textup{Aff}_k$
associates the largest Kan subcomplex of $\DP(X\times_kS)\simeq \Fun_k(\DP(X),\DP(S))$, where the equivalence is due to derived Morita equivalence.
Here $\Fun_k(\DP(X),\DP(S))$ denotes the $\infty$-category
of $\DP(k)$-linear exact functors. The stack $\mathcal{M}_{\textup{perf}}(X)$
is the (non-derived) moduli stack
of perfect complexes on $X$ (cf. \cite{TVM}).
The natural forgetful maps
\[
\Map^{\otimes}_k(\DP(X),\DP(S))\to \Map_k(\DP(X),\DP(S))
\]
(cf. Remark~\ref{monoidalqcperf}) induce a natural transformation
$\mathcal{F}_{X}\to \mathcal{M}_{\textup{perf}}(X)$, where we denote by $\Map_k(\DP(X),\DP(S))$ the largest Kan subcomplex of $\Fun_k(\DP(X),\DP(S))$.
Therefore the forgetful functor induces
a morphism
\[
X\longrightarrow \mathcal{M}_{\textup{perf}}(X).
\]
It is an analogue of the embedding $C\hookrightarrow \textup{Pic}(C)$
from an algebraic curve $C$ provided with a fixed point,
to its Jacobian.
\end{Remark}

Let $f:S\to \XX$ be a morphism of stacks.
Then we have an adjoint pair
\[
f^*:\DQ(\XX)\rightleftarrows \DQ(S):f_*.
\]
Conversely, when does an adjoint pair arise in such a way?
The following is a categorical characterization of functors
associated to morphisms $S\to \XX$, that is,
Theorem~\ref{supermain} implies a tannakian characterization theorem.

\begin{Theorem}
\label{geomono}
Let $\XX$ be an algebraic stack over $k$, that satisfies
either condition (i) or (ii) in Section \ref{subsec:2-3}.
Let $S$ be a scheme over $k$.
Let $\Phi:\DQ(\XX)\to \DQ(S)$ be a colimit-preserving
functor.
Then 
there exists a morphism $f:S\to \XX$ over $k$ such that
$f^*$ is equivalent to
$\Phi$
if and only if $\Phi$ is equivalent to a $k$-linear
symmetric monoidal functor
$($as objects in $\Map(\DQ(\XX),\DQ(S))$$)$.
\end{Theorem}

\begin{Corollary}
Let $\Psi:\DQ(S)\to \DQ(\XX)$ be a right adjoint functor i.e.,
an accessible and limit-preserving functor (see the
$\infty$-categorical adjoint functor theorem \cite[5.5.2.9]{HTT}).
Under the the same assumption as Theorem~\ref{geomono},
there is a $k$-morphism $f:S\to \XX$
such that $\Psi$ is equivalent to $f_*$ if and only if
a left adjoint $\Phi$ of $\Psi$ is equivalent to
the underlying functor of some $k$-linear symmetric monoidal functor
$\DQ^\otimes(\XX)\to \DQ^\otimes(S)$.
\end{Corollary}

\begin{Remark}
The above characterization
gives an answer to the question:
``what is the relationship between the group of automorphisms
of the derived category of a projective variety and the
group of isomorphisms of the variety?''
(see \cite[Preface]{BBH}).
An automorphism of a quasi-projective scheme $X$ over $k$
naturally corresponds to an autoequivalence of the symmetric monoidal
$\infty$-category $\DQ^{\otimes}(X)$ over $\DQ^\otimes(k)$.
It is perhaps worth remarking that Corollary~\ref{geomono} is new
even in the case $\XX$ is a scheme as well as the main theorem.
\end{Remark}

Acknowledgements.
We would like to experess our gratitude
to Yoshiyuki Kimura, Jacob Lurie, Shunji Moriya and Bertrand To\"en
for valuable comments.
Lurie informed us that he obtained a similar result in the setting
of derived schemes independently.
We also thank Hiroyuki Minamoto for stimulating conversations on this subject
and Toshiro Hiranouchi from whom we learned the works by Thomason and Balmer.
Finally,
we would like to experess our gratitude
to the referee for valuable comments.
The second author was supported by JSPS.


\begin{thebibliography}{99}





\bibitem{Bal}
P. Balmer,
Presheaves of triangulated categories and reconstruction of schemes,
Math. Ann., 324, (2002), 557-580.

\bibitem{Bal2}
P. Balmer,
The spectrum of prime ideals in tensor triangulated categories,
J. Reine Angew. Math., 588 (2005), 149--168.

\bibitem{BBH}
C. Bartocci, U. Bruzzo, and D. H. Ruip\'erez,
Fourier-Mukai and Nahm Transforms in Geometry and
Mathematical Physics,
Progress in Math. 279, (2009) Birkh\"auser.



\bibitem{BFN}
D. Ben-Zvi, J. Francis and D. Nadler,
Integral transforms and Drinfeld centers in derived algebraic geometry,
J. Amer. Math. Soc. (2010) 909--966.




\bibitem{Berg}
J. E. Bergner,
Three models for the homotopy theory of homotopy theories,
Topology 46
(2007), 397-436.


\bibitem{DM}
P. Deligne and J. Milne,
Tannakian categories,
Hodge Cycles, Motives, and Shimura Varieties (1982) 900. Berlin: Springer. 101--228.





\bibitem{Hop}
M. J. Hopkins,
Global methods in homotopy theory,
In Homotopy Theory (Durham, 1985),
London Math. Soc. Lecture Notes. vol.117 (1987), 73--96. 



\bibitem{Hov98}
M. Hovey,
Model Categories,
Mathematical Surveys and Monographs, vol. 63, American Mathematical Society, Providence, R.I., (1998).


\bibitem{Ill}
L. Illusie
Complex Cotangent et D\'eformations I
Lecture Notes in Mathematics 239, Springer-Verlag Berlin-Heidelberg-New York 1971.

\bibitem{Iwa}
I. Iwanari,
The category of toric stacks,
Compositio Math. 145, (2009) 718--746.

\bibitem{Tan}
I. Iwanari,
Tannakization in derived algebraic geometry,
preprint (2011) arXiv:1112.1761.





\bibitem{JJ}
A. Joyal,
Quasi-categories and Kan complexes,
J. Pure Appl. Algebra 175 (2002) 207--222.


\bibitem{J}
A. Joyal,
Notes on quasi-categories,
draft


\bibitem{JT}
A. Joyal and M. Tierney,
Quasi-categories vs Segal spaces,
Categories in algebra, geometry and mathematical physics, Contemp. Math., 431, Providence, R.I.: Amer. Math. Soc., 277--326

\bibitem{Kri}
A. Krishna,
Perfect complexes on Deligne-Mumford stacks and applications,
J. K-theory (2008) 1--45.

\bibitem{LM}
G. Laumon and L. Moret-Bailly,
Chapms alg\'ebriques,
Springer-Verlag (2000).


\bibitem{OL}
V. Lunts and D. Orlov,
Uniqueness of enhancement of triangulated categories,
J. Amer. Math. Soc. 23 (2010), 853--908.






\bibitem{HTT}
J. Lurie,
Higher Topos Theory,
Annals of Math. Studies, no. 170 (2009).

\bibitem{HA}
J. Lurie,
Higher Algebra,
preprint, February 2012 available at the author's webpage.

\bibitem{DAG8}
J. Lurie,
Derived algebraic geometry VIII,
preprint, (2011).

\bibitem{LT}
J. Lurie,
Tannaka duality for geometric stacks,
preprint.

\bibitem{DAG}
J. Lurie,
Derived algebraic geometry,
MIT Ph. D. Thesis.




\bibitem{Nee}
A. Neeman,
The chromatic tower for D(R), with an appendix by Marcel B\"okstedt,
Topology 31, (1992) 153--175.









\bibitem{Ols}
M. Olsson,
On proper coverings for Artin stacks,
Adv. Math. 198 (2005), no. 1, 93--106. 



\bibitem{Q}
D. Quillen,
Homotopical Algebra,
Lecture Notes in Math., 43 Springer-Verlag, Berlin and New York,
1967.


\bibitem{Saa}
N. Saavedra Rivano,
Categories Tannakiennes,
Lecture Notes in Mathematics 265, Springer-Verlag Berlin-Heidelberg-New York 1972.





\bibitem{Sav}
V. Savin,
Tannaka duality for quotient stacks,
manuscripta math., 119, (2006) 287--303.



\bibitem{Spal} 
N. Spaltenstein, 
Resolution of unbounded complexes. 
Compositio Math., 65, (1998) 121--154.




\bibitem{Tho}
R. Thomason,
The classification of triangulated subcategories,
Compositio Math., 105, (1997), 1--27.


\bibitem{TT}
R. Thomason, T. Trobaugh,
Higher algebraic K-theory of schemes and of derived categories,
The Grothendieck Festschrift, Vol. III, 
Progr. Math. 88, Birkhauser Boston, Boston, MA, (1990) 247--435.

\bibitem{To}
B. To\"{e}n,
The homotopy theories of dg-categories and derived Morita theory,
Invent. Math., 167 (2007), 615--667.


\bibitem{Azu}
B. To\"en,
Derived Azumaya's algebras and generators for twisted derived categories,
Invent . Math. to appear, arXiv:1002.2599.

\bibitem{TH}
B. To\"en,
Homotopical and Higher Categorical Structures in
Algebraic Geometry,
habilitation's thesis arXiv:math/0312262


\bibitem{TVM}
B. To\"en and M. Vaqui\'e,
Moduli of objects in dg-categories,
Ann. Sci. de l'ENS 40 (2007) 387--444.




\bibitem{TV0}
B. To\"{e}n and G. Vezzosi,
Homotopical algebraic geometry I: Topos theory,
Advances in Math. 193, (2005) 257--372.




\bibitem{TV}
B. To\"{e}n and G. Vezzosi,
Homotopical algebraic geometry II: geometric stacks and applications,
Memoirs of the AMS.


\bibitem{TVK}
B. To\"{e}n and G. Vezzosi,
A remark on K-theory and S-categories,
Topology 43, (2004), 765--791.


\bibitem{Tot}
B. Totaro,
The resolution property for schemes and stacks, J. Reine Angew. Math. 577 (2004), 1--22. 


\bibitem{Wall}
J. Wallbridge,
Tannaka duality over ring spectra,
preprint (2012) arXiv;1204.5787.



\bibitem{Yek}
A. Yekutieli,
Dualizing complexes, Morita equivalence and the derived Picard group of a ring,
J. London Math. Soc., 60 (1999), 723--746.









\end{thebibliography}
\end{document}